\documentclass[1p, preprint]{elsarticle}

\usepackage{amsmath}%
\usepackage{amsfonts}%
\usepackage{amssymb}%
\usepackage{mathrsfs}
\usepackage[capitalize]{cleveref}
\usepackage{overpic}
\usepackage{algpseudocode,algorithm}
\usepackage{url}

\newcommand{\R}{\mathbb{R}}
\renewcommand{\d}{\,\textup{d}}
\renewcommand{\L}{\mathcal{L}}
\newcommand{\F}{\bm{F}}

\newcommand{\f}{\bm{f}}
\newcommand{\G}{\bm{G}}
\renewcommand{\u}{\bm{u}}
\usepackage[utf8]{inputenc}
\usepackage{mathtools}
\usepackage{enumitem}

\usepackage{bm, booktabs, multirow}
\DeclareMathOperator*{\argmin}{arg\,min}
\DeclareMathOperator*{\argmax}{arg\,max}
\graphicspath{ {./images/} }

\usepackage{xcolor}
\newcommand{\edit}[1]{{\color{black} #1}} 

\journal{Computer Methods in Applied Mechanics and Engineering}

\begin{document}

\begin{frontmatter}

\title{Principled interpolation of Green's functions learned from data}

\author[CEE]{Harshwardhan Praveen\corref{mycorrespondingauthor}}
\cortext[mycorrespondingauthor]{Corresponding author}
\ead{hp477@cornell.edu}

\author[Cambridge]{Nicolas Boull\'e}
\ead{nb690@cam.ac.uk}

\author[CEE,CAM]{Christopher Earls}
\ead{earls@cornell.edu}

\address[CEE]{School of Civil and Environmental Engineering, Cornell University, Ithaca, NY, 14853, USA}
\address[Cambridge]{Isaac Newton Institute for Mathematical Sciences, University of Cambridge, Cambridge, CB3 0EH, UK}
\address[CAM]{Center for Applied Mathematics, Cornell University, Ithaca, NY, 14853, USA}

\begin{abstract}
We present a data-driven approach to mathematically model physical systems whose governing partial differential equations are unknown, by learning their associated Green's function. The subject systems are observed by collecting input-output pairs of system responses under excitations drawn from a Gaussian process. Two methods are proposed to learn the Green's function. In the first method, we use the proper orthogonal decomposition (POD) modes of the system as a surrogate for the eigenvectors of the Green's function, and subsequently fit the eigenvalues, using data. In the second, we employ a generalization of the randomized singular value decomposition (SVD) to operators, in order to construct a low-rank approximation to the Green's function. Then, we propose a manifold interpolation scheme, for use in an offline-online setting, where offline excitation-response data, taken at specific model parameter instances, are compressed into empirical eigenmodes. These eigenmodes are subsequently used within a manifold interpolation scheme, to uncover other suitable eigenmodes at unseen model parameters. The approximation and interpolation numerical techniques are demonstrated on several examples in one and two dimensions.
\end{abstract}


\end{frontmatter}

\section{Introduction}

It has been said that differential equations are the language of the universe \cite{Strogatz}. Indeed, most of our known physical laws are expressed mathematically, as rate forms \cite{Feynman}; thus these types of equations are of great practical interest. Also of interest is the study of symmetries and invariant structures that occur within the solution operators accompanying such governing differential equations \cite{Olver}. These structures appear as patterns that are observable within the \emph{Green's functions} associated with some system of interest \cite{Boulle2021}. Though very useful in forming fast forward solvers, and for affording mechanistic insight, Green's functions are not simple to find, in practice \cite{Evans2010}. In response to this difficulty, the current paper proposes an approach for learning \emph{empirical Green's functions} (EFGs) from observational data emanating from the response of some system of interest.

Previous investigators have learned solution operators, mapping model coefficient functions directly into solutions to governing equations, using convolutional neural networks and graph neural networks~\cite{Ying,anandkumar}. Some have even employed auto-encoder deep neural networks, as a kind of \emph{Koopman operator}, uncovering latent spaces where response data from weakly nonlinear systems may be lifted, so as to become somewhat linear; and thus amenable for consideration with a learned Green's function applied to this abstract, linear space~\cite{deepgreen}. The motivation for the foregoing important contributions to the literature appears to rest in either efficiently solving an inverse problem, or discovering an efficient to apply, forward model. In contrast to those contexts, others have employed deep rational neural networks \cite{rational} to uncover mechanistic insight from Green's functions learned from data \cite{Boulle2021}. There have also been other works on solving and learning partial differential equations \edit{with data-driven approaches (including neural networks)} \cite{raissi2018deep, raissi2018hidden, raissi2019physics, rudy2017data, schaeffer2017learning, berg2017neural, brunton2016discovering, Stephany2022, Bonneville2022}, but our work is focused on learning a Green's function.

\begin{figure}[htbp]
\centering
\begin{overpic}[width=\textwidth]{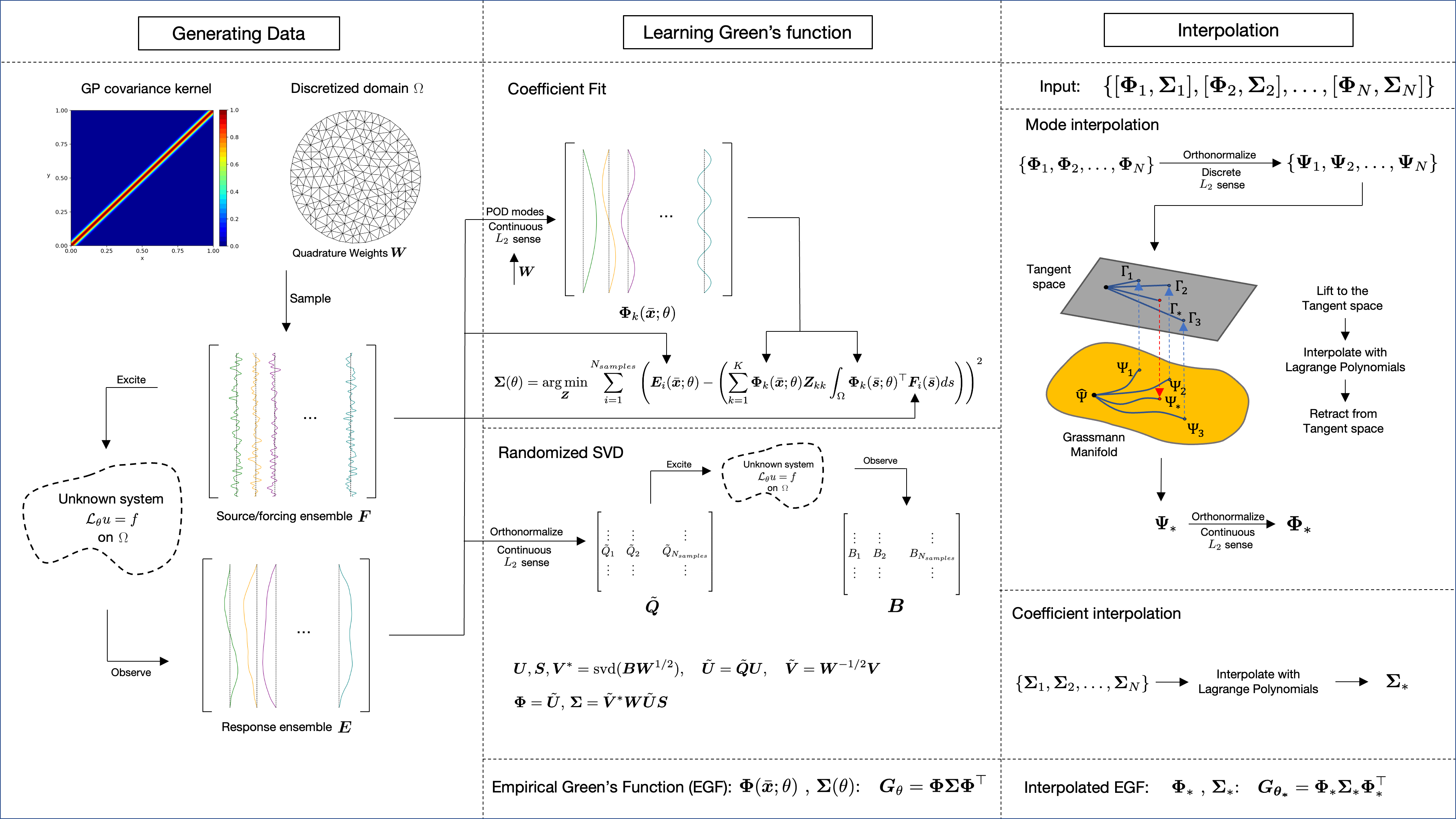}
\end{overpic}
\caption{Schematic of method for approximating and interpolating Green's functions associated with linear differential operators consisting of three steps: (1) Generation of the training data set using random functions sampled from a Gaussian process and associated solutions evaluated at sensor locations within the domain. (2) Construction of a low-rank approximant to the Green's function using POD-based coefficient fitting or the randomized singular value decomposition. (3) Interpolation between Green's functions to build an interpretable representation at unseen model parameter instances, $\theta\in\mathbb{R}$.}
\label{fig:method-schematic}
\end{figure}

In the present work, we propose a method for learning Green's functions from data without resorting to machine learning. This idea, illustrated in \cref{fig:method-schematic}, stems from earlier work aimed at learning Green's function-like, Fredholm integral kernels using environmental data, in order to formulate surrogate, reduced-order models exhibiting ``just enough physics'' to be useful in practical settings. The demonstration case for that earlier work was in predicting electro-magnetic ducting within the marine atmospheric boundary layer \cite{Earlsb,Earlsbb}. The currently proposed method builds from these earlier ideas, so as to now learn Green's function solution operators from observational data. While the method outlined herein is restricted to systems governed by self-adjoint differential operators, we note that many physical contexts are self-adjoint. Nonetheless, this is an important limitation to note, but we offer that our proposed approach for extending our method to weakly nonlinear (\emph{i.e.},~semi-linear) contexts \cite{Evans2010} may offset concern over these limitations, by extending applicability into the nonlinear regime. We do this by interpolating our learned Green's function models, in a structure preserving manner.

In the current work for uncovering meaningful descriptions of systems, in the form of solution operators to some self-adjoint physical context, we outline two data-driven approaches. The first one exploits the spectral decomposition of the kernel (Green's function) within a first-kind Fredholm integral equation, in order to substitute its eigenfunctions with proper orthogonal decomposition modes (empirical eigenfunctions), learned from observed system response data under a well characterized forcing/source term. The second method uses the randomized singular value decomposition (SVD) to construct a low-rank approximation for the Green's function, by using random Gaussian process with correlated entries, sampled from a multivariate normal distribution, to probe the unknown linear differential operator. We also propose a manifold interpolation scheme, for use in an offline-online setting, where offline excitation-response data, taken at specific model parameter instances, are compressed into empirical eigenmodes. These eigenmodes are subsequently used within a manifold interpolation scheme, to uncover other suitable eigenmodes, for an unseen model parameter instance; thus rendering an online, ``just-in-time'' EGF (obtained without the benefit of excitation-response data.) This interpolation approach is demonstrated on one and two-dimensional problems.

The present paper contains, in \cref{sec-methodology}, a complete discussion on the proposed methods for learning and interpolating empirical Green's Functions, along with details concerning considerations that are germane for the required data collection from the system of interest. Numerical examples of the method, applied to benchmark partial differential equations (PDEs), in one and two dimensions, with, and without noise appear in \cref{sec-numerical-examples}. Finally, concluding remarks are offered in \cref{sec-conclusion}.

\section{Methodology} \label{sec-methodology}

We consider a self-adjoint linear differential operator, $\mathcal{L}_{\bm{\theta}}$ depending on a set $\bm{\theta}=(\theta_1,\ldots,\theta_{N_{\text{params}}})\in\R^{N_{\text{params}}}$ of modelling parameters, defined on a bounded and connected domain $\Omega\subset \R^d$ in dimension $1\leq d\leq 3$, governing a physical system in the form of a boundary value problem:
\begin{equation} \label{problem}
\begin{aligned}
\mathcal{L}_{\bm{\theta}} u &= f, &&\text{on } \Omega, \\
\mathcal{B}(u) &= c, &&\text{in } \partial \Omega,
\end{aligned}
\end{equation}
where $\mathcal{B}$ is a linear differential operator specifying the boundary conditions of the problem, $f$ is a forcing (or source) term, and $u$ is the unknown system response. Under suitable conditions on the operator $\mathcal{L}_{\bm{\theta}}$ (\emph{e.g.},~uniformly elliptic), there exists a Green's function $G_{\bm{\theta}}:\Omega\times\Omega\to\R$ associated with $\mathcal{L}_{\bm{\theta}}$ that is the impulse response of the linear differential operator and defined \cite{Evans2010} as:~
\[\mathcal{L} G_{\bm{\theta}}(\vec{x},\vec{s}) = \delta (\vec{s} - \vec{x}), \quad \vec{x},\vec{s} \in \Omega,\]
where $\mathcal{L}$ acts on the first variable and $\delta$ is the Dirac delta function. Then with homogeneous Dirichlet boundary conditions, \emph{i.e.}, $\mathcal{B}(u) = 0$ on the boundary of the domain, the solution to \cref{problem}, for a forcing, $f$, can be expressed using the Green's function within a Fredholm integral equation of the first-kind:
\[u(\vec{x}) = \int_\Omega G_{\bm{\theta}}(\vec{x},\vec{s}) f(\vec{s}) d \vec{s}, \qquad \vec{x} \in \Omega.\]
Moreover, according to the theory of first-kind Fredholm integral equations, the integral operator associated with $G_{\bm{\theta}}(\vec{x},\vec{s})$ has the following spectral decomposition~\cite{Hsing}:
\[G_{\bm{\theta}}(\vec{x}, \vec{s}) = \sum_{k\geq 1} \frac{\Psi_k^* (\vec{x}; \bm{\theta}) \Psi_k (\vec{s}; \bm{\theta})}{\omega_k (\bm{\theta})},\quad \vec{x},\vec{s}\in\Omega,
\]
where $\omega_k(\bm{\theta})\in\mathbb{R}$ is the $k$th smallest eigenvalue (in absolute value) of $\mathcal{L}_{\bm{\theta}}$ associated with the eigenfunction, $\Psi_k$, and satisfying $\mathcal{L}_{\bm{\theta}}\Psi_k=\omega_k\Psi_k$ or, equivalently,
\[\int_{\Omega} G_{\bm{\theta}}(\vec{x},\vec{s}) \Psi(\vec{s};\bm{\theta}) \d \vec{s} = \frac{1}{\omega_k} \Psi_k(\vec{x};\bm{\theta}),\quad \vec{x}\in\Omega.\]
Then, the eigenfunctions of the Green's function coincide with the eigenfunctions of the self-adjoint operator $\mathcal{L}_{\bm{\theta}}$. We will now discuss a way to replace these unknown operator eigenfunctions with \textit{empirical eigenfunctions}; learned from the collected observational data from a system of interest. We begin with a description of the data (in our case synthetic data generated from simulations) which we use to learn the Green's functions.

\subsection{Generation of the training dataset} \label{subsec-building-ensemble}

The training datasets used to learn Green's functions consist of $N_{\text{samples}}\geq 1$ forcing terms, $\{f_j\}_{j=1}^{N_{\text{samples}}}$, and associated system's responses, $\{u_j\}_{j=1}^{N_{\text{samples}}}$, satisfying \cref{problem}. The selection of the forcing terms plays an important role in the accuracy of the algorithm that learns the Green's function. One would ideally want the set of forcing terms to be sufficiently ``diverse'' so that it approximates the vector space spanned by the dominant right singular vectors of the Green's functions~\cite{Boulle2021SVD}. For instance, if the forcing terms in the training dataset are orthogonal to the first right singular vector, then one cannot hope to approximate the largest singular value of the Green's function. To mitigate this risk, we probe the system~\eqref{problem} with random functions, sampled from a Gaussian process (GP), $\mathcal{GP}(0,K)$, with a covariance kernel, $K$, as motivated by recent theoretical results \cite{Boulle2021SVD,Boulle2021Theory}.

The forcing terms of the one-dimensional examples presented in \cref{sec-numerical-examples} are drawn from a GP, with mean zero and squared-exponential covariance kernel, $K_{\text{SE}}$, defined as
\[K_{\text{SE}}(x, y) = \exp \left(-|x-y|^2/(2 \sigma^2) \right), \quad x, y \in \Omega,\]
where $\Omega\subset\R$ is the problem domain and $\sigma>0$ is the length scale parameter. This parameter is chosen to be larger than the spatial discretization at which we are numerically evaluating the forcing terms, to ensure that the random functions are resolved properly within the discrete representation employed. \edit{Additionally, one should choose $\sigma$ sufficiently small so that the set of forcing terms is of full numerical rank.} \cref{sec-numerical-examples} will also feature numerical examples on a unit disk. In this case,  we generate the random forcing terms using the \texttt{randnfundisk} function of the Chebfun software system~\cite{chebfun,filip2019smooth} implemented in MATLAB. This function returns a smooth random forcing term defined on the unit disk, with a maximum frequency of approximately $2 \pi/\sigma$, and a standard normal distribution, $N(0,1)$, at each point on the disk. We note that the other types of covariance kernels employed in recent deep learning work~\cite{li2020fourier}, such as Green's functions associated with Helmholtz equations, are likely to yield better approximation results, because they already contain some information about the singular vectors of the operator, $\L_{\bm{\theta}}$. \edit{Hence, the theoretical bounds for the randomized SVD with arbitrary covariance kernels show that one may obtain higher accuracy by incorporating knowledge of the leading singular vectors of the differential operator into the covariance kernel~\cite{Boulle2021SVD}. In the present work, we employ a generic squared-exponential kernel as one may not have prior information about the governing operator, $\L_{\bm{\theta}}$, in real applications.}

For simplicity in notation, we now restrict to a scalar model parameter $\theta \in \R$, and an associated system response $u_j$, under a forcing term $f_j$; obtained by solving \cref{problem} using a finite element method with piecewise quadratic Lagrange polynomials implemented in the FEniCS software~\cite{logg2012automated}. Note that the methods we present here, including the manifold interpolation described in \cref{sec_manifold}, can be extended to work for operators instantiated with a collection of parameters, $\bm{\theta} = \{\theta_1, \theta_2, \dots, \theta_{N_{\text{param}}}\} \in \R^{N_{\text{param}}}$, but we now restrict ourselves to cases where $\theta \in \R$. We are interested in applications where one can only measure the responses at a finite number of locations, $\bar{\bm{x}} = \{\vec{x}_1, \ldots, \vec{x}_{N_{\text{sensors}}}\}$, where ${x}_i \in \Omega$ and $N_{\text{sensors}}\geq 1$ is the number of measurements taken within the domain, $\Omega$. For simplicity, the sensor locations in this work, $\{\vec{x}_i\}_{i=1}^{N_{\text{sensors}}}$, correspond with the nodes of the mesh that discretize the domain, but one could also evaluate the responses at arbitrary locations as well. We sample both the forcing terms and system responses at these points to assemble the following column vectors of dimension $N_{\text{sensors}}$:
\[
\f_i \coloneqq \begin{bmatrix} f_i(\vec{x}_1) & \cdots & f_i(\vec{x}_{N_{\text{sensors}}})\end{bmatrix}^\top, \quad
\u_i(\theta) \coloneqq \begin{bmatrix} u_i(\vec{x}_1) & \cdots & u_i(\vec{x}_{N_{\text{sensors}}})\end{bmatrix}^\top,
\]
The notation $\u_i(\theta)$ highlights the dependence of the responses on the parameter $\theta$, and $\top$ denotes the matrix transpose. These columns are then assembled to create sets of forcing terms, $\F(\bar{\bm{x}})$, and corresponding system responses, $\bm{E}(\bar{\bm{x}}, \theta)$, respectively, as

\[\F(\bar{\bm{x}}) \coloneqq \begin{bmatrix} \f_1 & \cdots & \f_{N_{\text{samples}}} \end{bmatrix}, \quad
\bm{E}(\bar{\bm{x}}; \theta) \coloneqq \begin{bmatrix} \u_1(\theta) & \cdots & \u_{N_{\text{samples}}}(\theta) \end{bmatrix},\]
where $\F(\bar{\bm{x}})$ and $\bm{E}(\bar{\bm{x}}; \theta)$ are real matrices of size $N_{\text{sensors}}\times N_\text{samples}$. The columns of the forcing matrix $\bm{F}$ are independent and identically distributed (\emph{i.i.d.}) and follow a multivariate normal distribution with covariance matrix $\bm{K}=(K(\vec{x}_i,\vec{x}_j))_{1\leq i,j\leq N_{\text{sensors}}}$.

\subsection{Low-rank approximation of Green's functions} \label{sec_low_rank_approx}

Let $\G_{\theta}=(G_{\theta} (\vec{x}_i,\vec{x}_j))_{1\leq i,j\leq N_{\text{sensors}}}$ be the symmetric matrix of the Green's function $G_{\theta}$, associated with the self-adjoint operator $\L_{\theta}$, and evaluated at the measurement locations $\bar{\bm{x}} = \{\vec{x}_1,\ldots,\vec{x}_{N_{\text{sensors}}}\}$. Given the sets $\F$ and $\bm{E}$ of forcing terms and responses, we are interested in computing a low-rank approximation to the Green's matrix $\G_{\theta}$, as
\begin{equation} \label{eq_low_rank_approx}
    \G_{\theta} \approx \bm{\Phi}\bm{\Sigma}\bm{\Phi}^\top,\quad
    \bm{\Phi} = \begin{bmatrix} \bm{\Phi}_1 & \cdots & \bm{\Phi}_K \end{bmatrix} \in \R^{N_{\text{sensors}}\times K},\quad \bm{\Sigma}=\textup{diag}(\sigma_1,\ldots,\sigma_K),
\end{equation}
where $K\geq 1$ is the target rank, and $|\sigma_1|\geq \cdots\geq |\sigma_K|$. We require our matrix of empirical eigenvectors $\bm{\Phi}$, to be orthonormal with respect to the quadrature weight matrix $\bm{W}=\textup{diag}((w_{i})_{i=i}^{N_{\text{sensors}}})$, associated with the finite element discretization of the forcing terms and solutions, \emph{i.e.}, $\bm{\Phi}^\top \bm{W}\bm{\Phi}=\bm{I}_K$, where $\bm{I}_K$ is the $K\times K$ identity matrix. We aim that the $k^{th}$ empirical eigenvector approximates the $k^{th}$ eigenfunction $\psi_k$ of the operator $\L_{\theta}$ in the following $L^2$-sense:
\[\lim_{N_{\text{sensors}} \to \infty} \sum_{i=1}^{N_{\text{sensors}}} w_i \left[\Psi_k(\vec{x}_i) - \bm{\Phi}_k(i) \right]^2= 0, \quad 1\leq k\leq K.\]
This ensures the $L^2$-convergence of the finite element function represented by the vector $\bm{\Phi}_k$ to $\Psi_k$, as the number of sensors increases. Once a low-rank approximant has been obtained, the discretized solution $u$ to \cref{problem}, under a new forcing term $f$, can be efficiently computed as
\begin{equation} \label{response-egf}
\bm{u} \approx \sum_{k=1}^K\sigma_k \bm{\Phi}_k \bm{\Phi}_k^\top \bm{W}\f, \quad \f=\begin{bmatrix}f(\vec{x}_1) & \cdots & f(\vec{x}_{N_{\text{sensors}}}) \end{bmatrix}^\top.
\end{equation}

\subsubsection{Proper orthogonal decomposition and least-square fitting of the coefficients} \label{subsec-fitting-green}

This section presents an algorithm for computing a low-rank approximation to $\bm{G}_\theta$, of the form of \cref{eq_low_rank_approx}, from pairs of forcing terms $\bm{F}$, and associate responses $\bm{E}$. The left singular vectors of the output ensemble $\bm{E}$ are used as the proper orthogonal decomposition (POD) modes that serve as our empirical eigenvectors $\bm{\Phi} \in \mathbb{R}^{N_{\text{sensors}} \times K}$, where $K$ is the desired rank of the approximant~ \cite{Liang2002}. The first step of the method consists of computing a singular value decomposition of the matrix $\bm{E}$, to obtain an orthonormal basis for the range of the integral operator associated with the Green's function:
\begin{align*}
\bm{U}, \bm{S}, \bm{V}^* &= \text{svd} (\bm{W}^{1/2} \bm{E}), \quad\bm{\Phi} = \bm{W}^{-1/2}\bm{U}_K,
\end{align*}
where $\bm{U}_K$ is the matrix obtained by truncating $\bm{U}$ to the first $K$ columns, and $\bm{\Phi}$ approximates the first $K$ eigenfunctions of the Green's function. We note that the empirical eigenvectors are orthonormalized using the quadrature weight matrix $\bm{W}$. Additionally, the matrix $\bm{S}$ of singular values provides us with a principled way to select the number of relevant empirical eigenvectors for constructing the \textit{empirical Green's function} (EGF) using the \emph{POD method}.

Once the empirical eigenvectors have been computed, one can follow \cref{eq_low_rank_approx} to obtain a rank $K$ empirical Green's function $\bm{G}_{\theta}$ for the specified sensor locations, as $\bm{G}_{\theta}\approx \bm{\Phi}\bm{\Sigma}\bm{\Phi}^\top$, where $\bm{\Sigma}\in \mathbb{R}^{K\times K}$ is an unknown diagonal coefficient matrix. We approximate the coefficients of the diagonal matrix $\bm{\Sigma}(\theta)$ by solving the following least squares problem:
\begin{equation} \label{minimization-eq}
\bm{\Sigma} \approx \argmin_{\bm{Z}=\textup{diag}((Z_k)_{k=1}^K)} \sum_{i=1}^{N_{\text{samples}}}\left\|\bm{u}_i-\bm{\Phi}\bm{Z}\bm{\Phi}^\top\bm{W}\bm{f}_i\right\|_2^2,
\end{equation}
where $\|\cdot\|_2$ denotes the matrix $2$-norm. We construct the normal equations for the least square minimization and then solve the system of equations using the linear algebra library in NumPy~\cite{harris2020array}. The main advantage of this algorithm is that it only requires a single pass over the differential operator, and does not assume any distribution on the forcing terms. In the next section, we will introduce a second algorithm which allows us to obtain higher accuracy under the assumption that we have control over the forcing terms that act on the system.

\subsubsection{Randomized singular value decomposition} \label{subsec-randomized-svd}

The randomized singular value decomposition (SVD) is a popular algorithm for computing a low-rank approximant to a large matrix $\bm{G}\in\R^{N_{\text{sensors}}\times N_{\text{sensors}}}$ using matrix-vector products with random vectors $\bm{f}_{1},\ldots,\bm{f}_{N_{\text{samples}}}\in \R^{N_{\text{sensors}}}$~\cite{halko2011finding,martinsson2020randomized}. The error analysis for the randomized SVD in~\cite{halko2011finding} uses standard Gaussian random vectors, but other random embedding techniques have been considered such as random permutations~\cite{ailon2009fast}, sparse sign matrices~\cite{clarkson2017low,meng2013low,Nelson2013OSNAP,urano2013fast}, and subsampled randomized trigonometric transforms (SRTTs)~\cite{ailon2009fast,ailon2006approximate,Parker95randombutterfly,woolfe2008fast}, to mitigate the computational cost of Gaussian vectors in large-scale numerical linear algebra applications. Here, we employ a generalization of the randomized SVD, which uses random Gaussian vectors with correlated entries, sampled from a multivariate normal distribution~\cite{Boulle2021SVD}. Hence, we are interested in probing an unknown linear differential operator using smooth random forcing terms sampled from a Gaussian process, which implies that the discretized forcing term $\bm{F}$ has correlated rows but independent columns (see \cref{subsec-building-ensemble}).

We mainly use the algorithm described in \cite[Sec.~1.5]{halko2011finding} with a slight modification, to ensure that the approximated singular modes are orthonormal with respect to the weight matrix $\bm{W}$ associated with the finite element discretization. The first step in the \emph{randomized SVD method} consists of computing the economized QR decomposition of the system's response matrix $\bm{E}$, to obtain an orthonormal basis for the range of $\bm{G}$:
\[\bm{Q},\bm{R}=\text{qr}(\bm{W}^{1/2}\bm{E}),\quad \tilde{\bm{Q}}=\bm{W}^{-1/2}\bm{Q}.\]
The multiplication by the weight matrix ensures that the columns $\tilde{q}_1,\ldots,\tilde{q}_{K}$ of $\tilde{\bm{Q}}$ are orthonormal in the $L^2$-norm associated with the finite element discretization, \emph{i.e.},
\[\sum_{j=1}^{N_{\text{sensors}}} w_j \tilde{q}_{k}(j)\tilde{q}_{k'}(j) = \delta_{kk'},\quad 1\leq k,k'\leq K,\]
where $\delta_{kk'}$ is the Kronecker delta symbol. The next step consists of forming the matrix $\bm{B}=\tilde{\bm{Q}}^*\bm{G}=(\bm{G}^*\tilde{\bm{Q}})^*\in\R^{K\times N_{\text{sensors}}}$, which requires the solution of the adjoint problem to \cref{problem}, under forcing that is prescribed by the columns of $\tilde{\bm{Q}}$. In this work, we assume that the partial differential operator is self-adjoint, and solve \cref{problem} with the forcing terms $\tilde{q}_1,\ldots,\tilde{q}_{K}$. We then evaluate the corresponding solutions at the nodes of the mesh to form $\bm{B}$. Finally, we compute the singular value decomposition of the matrix $\bm{B}$ before computing the left and right singular vectors, $\tilde{\bm{U}}$ and $\tilde{\bm{V}}$, respectively, which approximate the singular vectors of $\tilde{\bm{G}}$ in the $L^2$ norm, as
\[\bm{U},\bm{S},\bm{V}^*=\text{svd}(\bm{B}\bm{W}^{1/2}),\quad \tilde{\bm{U}}=\tilde{\bm{Q}}\bm{U},\quad \tilde{\bm{V}}=\bm{W}^{-1/2}\bm{V}.\]
Since the partial differential operator is self-adjoint, we select the empirical eigenvectors to be the columns of $\tilde{\bm{U}}$, \emph{i.e.},~$\bm{\Phi} = \tilde{\bm{U}}$. We then choose the empirical eigenvalue matrix $\bm{\Sigma}\in\R^{K\times K}$ such that $\tilde{\bm{U}}\bm{S}\tilde{\bm{V}}^*=\tilde{\bm{V}}\bm{\Sigma}\tilde{\bm{V}}^*$, \emph{i.e.}, $\bm{\Sigma}=\tilde{\bm{V}}^*\bm{W}\tilde{\bm{U}}\bm{S}$, to account for the sign flip between the left and right singular vectors associated with the operator $\L_{\theta}$.

The main difference between this algorithm and the POD method described in \cref{subsec-fitting-green}, is that the randomized SVD requires two passes to probe the system: a first pass with random functions sampled from a Gaussian process, and a second pass with the functions defined by the columns of $\tilde{Q}$. While the randomized SVD has near-optimal theoretical guarantees and requires many fewer response samples than the POD algorithm, one may not be able to employ forcing terms associated with the columns of $\tilde{Q}$ in practical applications. However, several single pass randomized algorithms~\cite{tropp2017practical,tropp2019streaming,upadhyay2016fast} have been proposed to alleviate this issue and only require Gaussian input vectors, such as the generalized Nystr{\"o}m method~\cite{nakatsukasa2020fast,nystrom1930praktische,williams2000using}.

We remark that the approach described in this section could be generalized to non-self-adjoint operators by automatically deriving the adjoint associated with the partial differential operator $\L_{\theta}$ using the dolfin-adjoint software package~\cite{farrell2013automated}, and then solving the adjoint equation to compute an approximant to the Green's function. However, we are mainly motivated by applications where the operator is not known, and do not consider the non self-adjoint case in this paper. Finally, one may also exploit the hierarchical low-rank structure of partial differential operators~\cite{bebendorf2003existence,bebendorf2008hierarchical} to build a good approximation to the Green's functions with fewer forcing terms~\cite{Boulle2021Theory,boulle2022learning,lin2011fast,martinsson2019fast}.

\subsection{Interpolation of EGFs to unobserved model parameter instances} \label{sec-interpolation}

Since Green's functions assume linearity of the underlying differential operator governing the system, we extend our method to ``weakly nonlinear'' or semi-linear~\cite{Evans2010} contexts by assuming that the system responses at some model parameter instance, $\theta$, are locally linear, and thus suitable for forming an EGF that applies locally within an underlying nonlinear manifold. Again, we would like to emphasize that this interpolation can be extended to contexts where the models are instantiated by a set of parameters, $\bm{\theta}$. We navigate this space by moving from linear ``patch'' to linear ``patch'' using manifold interpolation. \edit{The interpolation method described in this section is inspired by~\cite{Amsallem:2008}, and extended due to~\cite{Absil}. We describe an algorithm (see \cref{interpolation-scheme}) which takes a set of $N_{\text{models}}$ eigenbases and coefficients, $\{\bm{\Phi}(\bar{\bm{x}};\theta_j), \bm{\Sigma}(\theta_j),\theta_j\}_{j=1}^{N_{\text{models}}}$, and interpolates between them to approximate a model at the target parameter instance $\theta_*$, in the form of interpolated eigenmodes, $\bm{\Phi}(\bar{\bm{x}};\theta_*)$, and corresponding coefficients, $\bm{\Sigma}(\theta_*)$.}

\begin{algorithm}[htbp]
\caption{Interpolation of EGFs to unseen modeled parameters}
\label{interpolation-scheme}
\begin{algorithmic}
\State \textsl{Input:} A set of eigenbases and coefficients, $\{\bm{\Phi}(\bar{\bm{x}};\theta_j), \bm{\Sigma}(\theta_j),\theta_j\}_{j=1}^{N_{\text{models}}}$ to be interpolated for a new parameter instance, $\theta_*$.
\State \textsl{Step 1:} From among the $\theta_j$'s, identify a $\hat{\theta}$, as being closest to $\theta_*$, and use the associated eigenbasis, $\bm{\Phi}(\bar{\bm{x}};\hat{\theta})$, as the origin point for the interpolation.
\State \textsl{Step 2:} Correct sign flips and shuffling of the eigenmodes, following the procedure described in \cref{sec_correct_sign_order}.
\State \textsl{Step 3:} Construct $L^2$-orthonormal matrices of eigenmodes using the quadrature weight matrix: $\bm{\Psi}(\bar{\bm{x}};\theta_j) = \bm{W}^{1/2} \bm{\Phi}(\bar{\bm{x}};\theta_j)$.
\State \textsl{Step 4:} Perform the interpolation:
\begin{enumerate}[noitemsep,nolistsep]
\item Lift to the tangent space $\mathcal{T}_{\hat{\bm{\Psi}}} \mathcal{G}_{n,r}$ (bijective to the horizontal space of $\mathcal{T}_{\hat{\bm{\Psi}}} \mathcal{ST}_{n,r}$) by using the following map \cite{Absil},
\[\bm{\Gamma}_j = \bm{\Psi}(\bar{\bm{x}};\theta_j) - \bm{\Psi}(\bar{\bm{x}};\hat{\theta})\ \text{sym}\left(\bm{\Psi}^\top(\bar{\bm{x}};\hat{\theta}) \bm{\Psi}(\bar{\bm{x}};\theta_j)\right), \quad \text{sym}\left(\bm{Y}\right) \coloneqq (\bm{Y} + \bm{Y}^\top)/2.\]
\item Using Lagrange polynomials, compute the interpolated tangent vector, $\bm{\Gamma}_*$, using $ \{\bm{\Gamma}_j\}_{j = 1}^{N_{models}}$. Interpolate the coefficients, $ \{\tilde{\bm{\Sigma}}(\theta_j)\}_{j = 1}^{N_{models}}$, with the same scheme, to obtain $\bm{\Sigma}(\theta_*)$.
\item Compute the interpolated eigenvectors, $\bm{\Psi}(\bar{\bm{x}};\theta_*)$, by mapping $\bm{\Gamma}_*$ back to $\mathcal{ST}_{n,r}$ using the exponential map \cite{Absil},
\[\bm{\Psi}(\bar{\bm{x}};\theta_*) = \text{qf}(\bm{\Psi}(\bar{\bm{x}};\hat{\theta}) + \bm{\Gamma}_*),\]
where $\text{qf}(\bm{A})$ denotes the $\bm{Q}$ factor of the QR decomposition of $\bm{A} \in \mathbb{R}^{n \times r}$.
\end{enumerate}
\State \textsl{Step 5:} $L^2$-orthonormalize $\bm{\Psi}$, using $\bm{W}$: $\bm{\Phi}(\bar{\bm{x}};\theta_*) = \bm{W}^{-1/2}\bm{\Psi}(\bar{\bm{x}};\theta_*)$.
\State \textsl{Step 6:} Order the eigenbasis and coefficients to match the eigenbasis at the initial parameter instance, $\hat{\theta}$.
\State \textsl{Output:} Interpolated eigenmodes $\bm{\Phi}(\bar{\bm{x}};\theta_*)$ and coefficients $\bm{\Sigma}(\theta_*)$ at parameter $\theta_*$.
\end{algorithmic}
\end{algorithm}

\subsubsection{Correcting eigenmodes sign and order} \label{sec_correct_sign_order}

The approximated Green's function (EGF), $\bm{G_\theta} (\bar{\bm{x}}, \bar{\bm{s}}) = \bm{\Phi}(\bar{\bm{x}}; \theta) \bm{\Sigma}(\theta) \bm{\Phi}(\bar{\bm{x}}; \theta)^\top$, constructed in \cref{sec_low_rank_approx}, is not affected by sign flips of its eigenmodes, as the normalized eigenvalue decomposition is unique up to a sign of the eigenvectors. Therefore, two Green's functions approximated at close parameters, $\theta_1$ and $\theta_2$, may have learned eigenmodes with opposite signs. Before performing the interpolation of the EGFs, we account for any potential sign flips by matching the signs of the individual eigenmodes. We first consider the eigenvectors at the origin of the manifold, $\bm{\Phi}(\bar{\bm{x}};\hat{\theta})$, around which we form the tangent space. We then compute an inner product between the eigenvectors of the other interpolants $\bm{\Phi}_k(\bar{\bm{x}};\theta_j)$ and the corresponding eigenvectors at the origin, $\bm{\Phi}_k(\bar{\bm{x}};\hat{\theta})$, in order to update the signs of the interpolants' eigenmodes, as follows:
\[\bm{\Phi}_k(\bar{\bm{x}};\theta_j) =
\begin{cases}
-\bm{\Phi}_k(\bar{\bm{x}};\theta_j), & \text{if}\ \langle \bm{\Phi}_k(\bar{\bm{x}};\theta_j), \bm{\Phi}_k(\bar{\bm{x}};\hat{\theta}) \rangle < 0,\\
\bm{\Phi}_k(\bar{\bm{x}};\theta_j), &\text{otherwise},
\end{cases}\]
where $\langle\cdot,\cdot\rangle$ denotes the inner product with respect to the quadrature weight matrix, $\bm{W}$.

\begin{figure}[htbp]
\centering
\begin{overpic}[width=\textwidth]{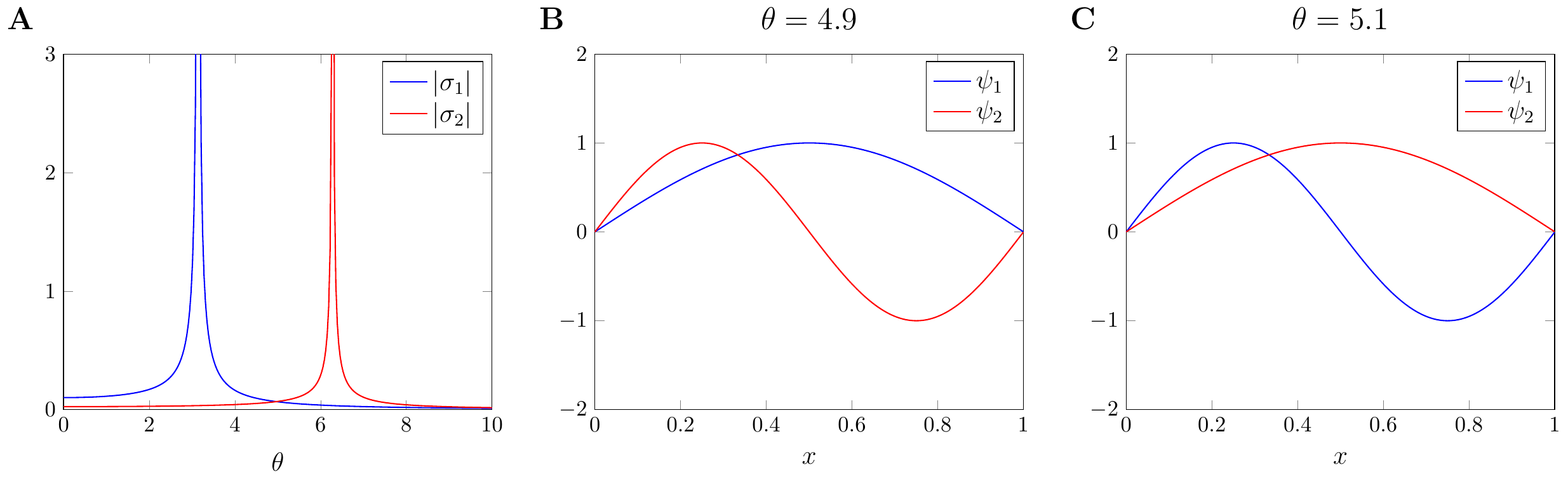}
\end{overpic}
\caption{(A) Magnitude of the first two eigenvalues of the Green's function associated with the one-dimensional Helmholtz equation. The first two corresponding eigenmodes are swapped when the frequency, $\theta$, increases beyond the critical value of $\theta_{\text{crit}} = \sqrt{5/2}\pi\approx 5$, as displayed in (B) and (C).}
\label{fig:mode-swap}
\end{figure}

Another potential issue that might arise during the manifold interpolation is the ``shuffling'' of eigenmodes. By construction, the eigenmodes approximated by the algorithms described in \cref{sec_low_rank_approx} are ordered according to the magnitude of the eigenvalues, but swapping of the eigenmodes may occur as the parameter $\theta$ varies. As an example, consider the one-dimensional Helmholtz equation with frequency $\theta\geq 0$ and homogeneous boundary conditions:
\[\frac{d^2u}{dx^2}+\theta^2 u=f,\quad x\in[0,1].\]
The first two eigenvalues of the associated Green's functions are given by $\sigma_1 = 1/(\theta^2-\pi^2)$ and $\sigma_2 = 1/(\theta^2-4\pi^2)$, and satisfy $|\sigma_1|>|\sigma_2|$, when $\theta<\sqrt{5/2}\pi$, and $|\sigma_1|\leq|\sigma_2|$ otherwise. This implies that the first two eigenmodes are swapped when $\theta>\sqrt{5/2}\pi$, as illustrated by \cref{fig:mode-swap}, and that the manifold interpolation technique will perform poorly in such cases. Note that the same phenomena occurs for the higher eigenmodes at larger values of $\theta$. We resolve this issue by reordering the eigenmodes of the discovered basis at the given interpolation parameter to match the ones at the origin point $\hat{\theta}$, where we are linearizing our Green's function. \edit{For a given parameter, $\theta_j$, and mode number, $1\leq k\leq K$, we select the $k$th eigenmode at $\theta_j$ to be the one with minimal angle with the $k$-th eigenmode, at parameter $\hat{\theta}$, \emph{i.e.},
\begin{equation*}
    \begin{aligned}
        &\bm{\Phi}_k(\bar{\bm{x}};\theta_j) = \bm{\Phi}_{\ell'}(\bar{\bm{x}};\theta_j),    \quad \text{where}\quad \ell' = \text{match}(k), \\
        \text{match}(k) = &\argmax_{1\leq \ell\leq K, \ell \ne \text{match}(k') \forall 1 \le k' < k} |\langle \bm{\Phi}_\ell(\bar{\bm{x}};\theta_j), \bm{\Phi}_k(\bar{\bm{x}};\hat{\theta}) \rangle|,\quad 1\leq k\leq K.
    \end{aligned}
\end{equation*}
Note that the reordering is done by matching the eigenmodes, starting from the eigenmode with the highest singular value, without replacement. When the order of the eigenmodes has been changed by this procedure, we also re-order the corresponding eigenvalues to preserve the value of the Green's function.}

\subsubsection{Manifold interpolation} \label{sec_manifold}

We begin by summarizing some useful notions and results from differential geometry; the interested reader may refer to~\cite{Absil} for a more exhaustive treatment. The Grassmann manifold $\mathcal{G}_{n,r}$ is defined as the set of all $r$-dimensional subspaces in $ \mathbb{R}^{n}$, and a particular $r$-dimensional subspace from $ \mathbb{R}^{n} $ is represented with an element from $\mathcal{G}_{n, r}$. Such an element can also be non-uniquely represented by some particular matrix, having orthonormal columns, $\bm{\Psi} \in \mathbb{R}^{n \times r}$: itself an element from within the equivalence class of matrices spanning the $r$-dimensional subspace in question. Such a matrix is an element from the \emph{compact Stiefel manifold}: the set of all such orthonormal matrices within $\mathbb{R}^{n \times r}$~\cite{Absil}. In our case, we construct each orthonormal matrix $\bm{\Psi}$ from data using a set of empirical eigenmodes $\bm{\Phi}$. Then, $\bm{\Psi}$ is an element within the compact Stiefel manifold $\mathcal{S}_{n,r}$. Our aim is to use pre-existing offline observations to create a collection of $\bm{\Phi}$'s that are suitable for an online interpolation, at a point where we have no observational data. It is not possible to directly interpolate the empirical eigenmodes modes within the compact Stiefel manifold since it is not a vector space. However, at each manifold point there exists a tangent space, which is a vector space, with its \emph{origin} occurring at the point of tangency. We use the appropriate tangent space to interpolate the offline empirical eigenmodes, to create a ``just-in-time,'' online mode set, for use in constructing a suitable Green's function at some desired model parameter instance $\theta_*$.

\begin{figure}[htbp]
\centering
\begin{overpic}[width=0.5\textwidth]{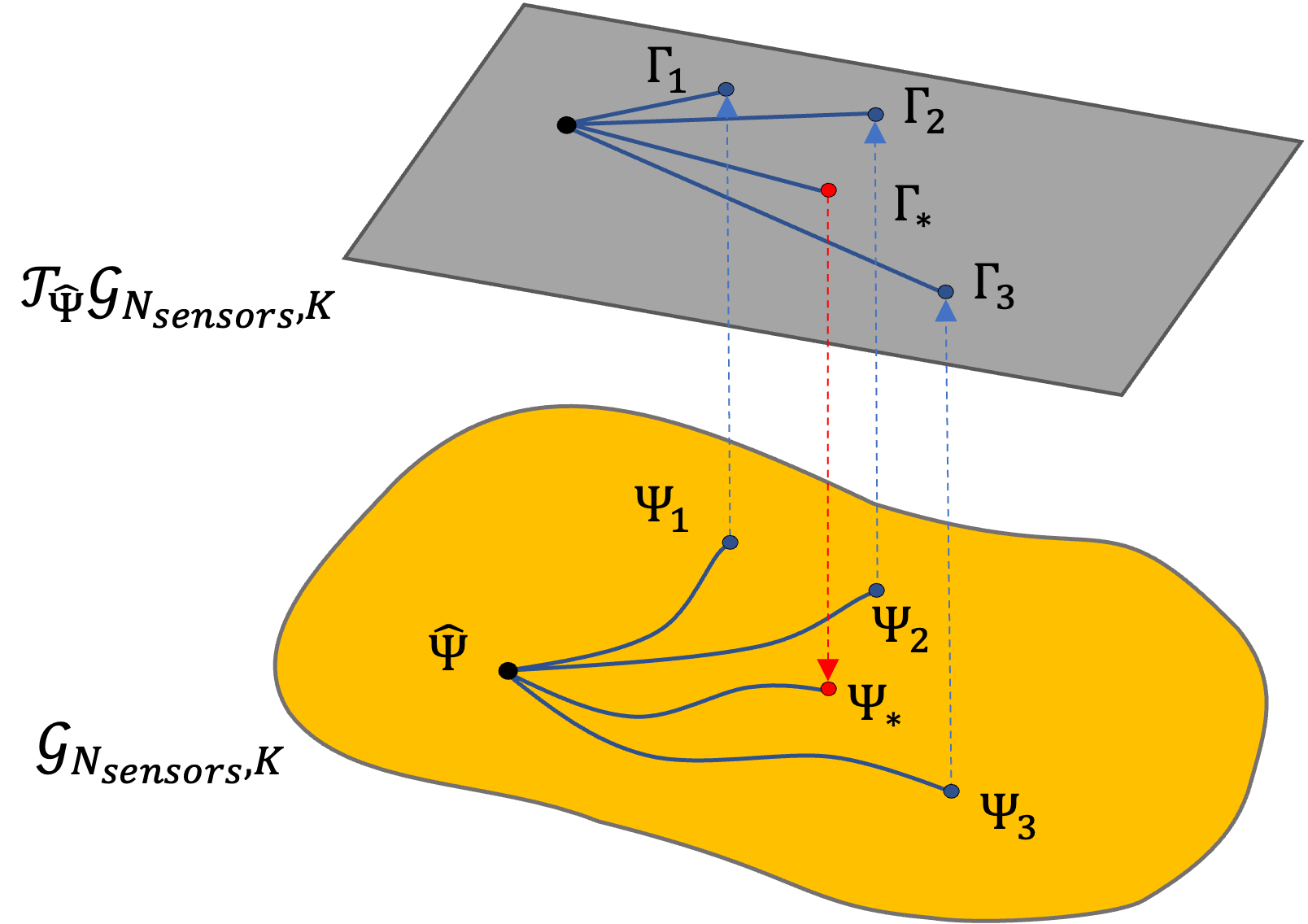}
\end{overpic}
\caption{The points on the Grassmann manifold (represented by the orthonormal matrices $\bm{\Psi}_j$) are projected onto the flat tangent space $\smash{\mathcal{T}_{\hat{\bm{\Psi}}} \mathcal{G}_{N_{\text{sensors}},K}}$ at the ``origin'' $\smash{\hat{\bm{\Psi}} = \bm{\Psi}(\bar{\bm{x}};\hat{\theta})}$, where they are interpolated and then returned to the Grassmann manifold as $\bm{\Psi}_*$.}
\label{fig:Interp}
\end{figure}

The manifold interpolation scheme we employ exploits the bijective relation between the tangent space to the Grassmann manifold, $\mathcal{T}_{\bm{\Psi}} \mathcal{G}_{n,r}$, and the horizontal space, within the tangent space to the compact Stiefel manifold, $\mathcal{T}_{\bm{\Psi}} \mathcal{S}_{n,r}$. Along with the preserved metric structure between these two spaces~\cite{Absil}, the use of the horizontal space has more practical benefits: it ensures uniqueness within the computational framework underpinning the interpolation, by allowing a particular matrix representation, from within the previously mentioned equivalence class, to be used in the computations~\cite{Edelman1998}. \cref{fig:Interp} offers a schematic overview of the manifold interpolation, while \cref{interpolation-scheme} describes the QR-decomposition based variant of the interpolation algorithm implemented in this work~\cite{Earlsb,Earlsa}. The scheme used in the present work relies only on a single QR decomposition. Thus, it is computationally cheaper than previous approaches~\cite{Amsallem:2008}, which require multiple singular value decompositions and matrix inversions. The process of returning from the horizontal space of $\mathcal{T}_{\bm{\Psi}} \mathcal{S}_{n,r}$, back to $\mathcal{S}_{n,r}$, with schemes that employ singular value decomposition might lead to a disruption within the ordering of the columns within the resulting matrix. For the problems that we consider, the resulting interpolated eigenbasis. using the proposed QR-based maps, are free from this mode ``shuffling''; though we still check for this for the reasons mentioned in \cref{sec_correct_sign_order}.

We remark that the matrix of empirical eigenmodes $\bm{\Phi}(\bar{\bm{x}};\theta)$ constructed in \cref{subsec-fitting-green,subsec-randomized-svd}, is orthonormal with respect to the quadrature weights, as $\bm{\Psi}(\bar{\bm{x}};\theta) = \bm{W}^{1/2} \bm{\Phi}(\bar{\bm{x}};\theta)$. We also need to specify the point on $\mathcal{T}_{\bm{\Psi}} \mathcal{S}_{n,r}$, where we construct the required tangent space, within whose horizontal space we perform the interpolation. For this, we select the mode set with the parametric value, $\theta_j$, closest to the target parameter point $\theta_*$, \emph{i.e.},
\[\hat{\theta} = \argmin_{1 \le j \le N_{models}} |\theta_j - \theta_*|,\]
where $\hat{\theta}$ denotes the parametric value at the origin. Once the origin has been identified, we use linear Lagrange polynomials to interpolate the tangent vectors $\{\bm{\Gamma}_j\}_{j = 1}^{N_{\text{models}}}$ within the horizontal space, as described in \cref{interpolation-scheme}. Since the associated coefficient matrices $\bm{\Sigma}(\theta_*)$ occur within a Euclidean space, manifold interpolation is not required. Therefore, we interpolate these directly with linear Lagrange polynomials.

\edit{As a closing remark, we note that in some cases it is possible to interpolate the matrices $\bm{\Phi}$, which correspond to the eigenmodes in our Green's function model, directly with simple linear interpolation (\emph{e.g.}, with splines). This considerably simplifies the procedure outlined in \cref{interpolation-scheme}. However, such linear interpolation, unlike the manifold-based interpolation described in this work, is not structure preserving~\cite{Degroote2010}, and its accuracy is very problem dependent.}

\section{Numerical results} \label{sec-numerical-examples}

In this section, we evaluate the algorithms described in \cref{sec-methodology}, for approximating Green's functions (as EGFs) from input-output pairs, and also interpolating to unseen parameter points, $\theta$. A number of synthetic problems, in one and two dimensions, are considered. Unless specified otherwise, we generate the different training data sets presented in this section using the parameters $N_{\text{sensors}}=2000$, $\sigma=5\times 10^{-3}$ (\emph{i.e.},~length-scale for the GP squared exponential covariance kernel), $N_{\text{samples}}=2000$, and apply the low-rank approximation algorithms with a target rank of $K=100$. Additionally, we only use $100$ training pairs in the randomized SVD case to illustrate its superiority for approximating Green's function with smaller training sets (at the cost of using a two pass procedure). For problems where the closed-form Green's function is not known, we estimate the relative error of the EGF empirically by generating a testing data set of $N_{\text{test}}=100$ input-output pairs of the form $\smash{\{(f_j,u_j)\}_{j=1}^{N_{\text{test}}}}$. Then, we define the testing error as
\[\epsilon_{\text{test}}=\frac{1}{N_{\text{test}}}\sum_{j=1}^{N_{\text{test}}}\frac{\|u_j-\tilde{u}_j\|_{L^2(\Omega)}}{\|u_j\|_{L^2(\Omega)}},\quad \tilde{u}_j=\int_{\Omega}G(\vec{x},\vec{s})f_j(\vec{s})\d \vec{s},\]
where the integrals are approximated by a quadrature rule that uses the finite element quadrature weights. Here, $\tilde{u}_j$ denotes the reconstructed solution from the learned Green's function. Note that we often express the relative error as a percentage by multiplying it by one hundred.

\subsection{Approximation of Green's functions}

We first evaluate the algorithms for approximating Green's functions from a training dataset of forcing terms and solutions on one and two-dimensional problems.

\subsubsection{One dimensional Poisson equation} \label{sec_1d_poisson}

We begin with a one-dimensional Poisson equation with homogeneous boundary conditions:
\[ -\frac{d^2u}{dx^2} = f, \quad x \in [0,1],\quad u(0)=u(1)=0.\]
The Green's function associated with this problem is available in closed-form as
\[G_{\text{exact}}(x,y) =
\begin{cases}
x (1 - s), & \text{if}\ x \le s, \\
s (1 - x), & \text{if}\ s > x,
\end{cases}\]
where $x, s \in [0,1]$.

\begin{figure}[htbp]
\centering
\begin{overpic}[width=0.8\textwidth]{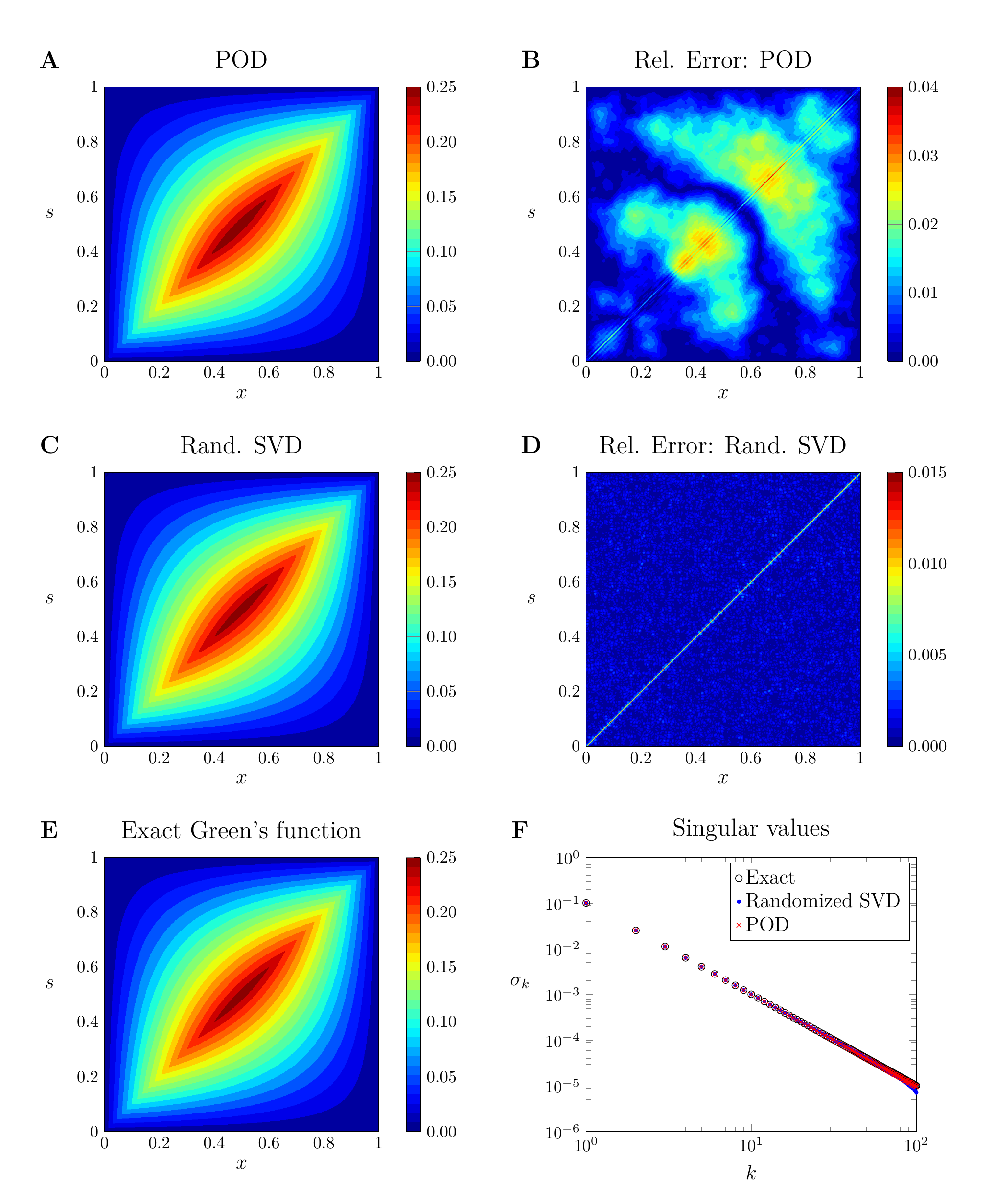}
\end{overpic}
\caption{EGF models for the 1D Poisson problem. (A) Green's function learned using the POD method. (B) Error contour for the POD method. (C) Green's function learned using the randomized SVD. (D) Error contour for the randomized SVD. (E) Exact Green's function. (F) Comparison of coefficients learned by both methods with the exact singular values.}
\label{poisson1D}
\end{figure}

In \cref{poisson1D}, we display the exact Green's function (panel E) along with its approximations constructed using the POD method (A) and the randomized SVD (C), introduced in \cref{subsec-fitting-green,subsec-randomized-svd}. We first compute the relative error $\epsilon$ between the learned EGFs $\G$, and the exact Green's function $\G_{\text{exact}} = G_{\text{exact}}(\bar{\bm{x}},\bar{\bm{s}})$ (evaluated at the sensor locations used to learn the EGF), as
\[\epsilon=\frac{\|\G-\G_{\text{exact}}\|_{L^2([0,1]\times[0,1])}}{\|\G_{\text{exact}}\|_{L^2([0,1]\times[0,1])}}\times 100,\]
and find that the POD method achieves a relative error of $1.1\%$, while the randomized SVD yields an error of $0.09\%$. We note that the performance of the randomized SVD is consistent with the theory~\cite{Boulle2021SVD,halko2011finding}, which predicts a relative error of the order of $0.01\%$ for a target rank of $K=100$ and Green's functions singular values, $\sigma_k=1/(\pi^2 k^2)$. We observe that the randomized SVD achieves a lower relative error than the POD method, while using fewer input-output training pairs. We also present the error contours for the approximated Green's function where we plot the relative error,
\[\epsilon_{\G}=\frac{|\G-\G_{\text{exact}}|}{\|\G_{\text{exact}}\|_{L^2([0,1]\times[0,1])}},\]
for the POD method in \cref{poisson1D}(B) and the randomized SVD (D). Lastly, panel (F) compares the learned singular values with the first hundred exact singular values associated with the Green's functions. While the randomized SVD yields a lower relative error on the Green's function, the POD method is surprisingly able to recover the higher order singular values more accurately, while also capturing the algebraic decay rate. Finally, we evaluate the relative error of the learned Green's function on the testing dataset, and obtain a testing error of $1.6\%$ for the POD method and $0.1\%$ for the randomized SVD.

\subsubsection{Noisy Poisson equation} \label{poisson1D-noisy-section}

We evaluate the effect of noise on the accuracy of the learned EGFs by perturbing the training solutions to the Poisson equation, generated in \cref{sec_1d_poisson}, with $10\%$ Gaussian white noise. We assume that our system only has access to forcing terms $\{f_j\}_{j=1}^{N_{\text{samples}}}$ and associated noisy solutions $\{u_j^{\text{noisy}}\}_{j=1}^{N_{\text{samples}}}$ satisfying
\[u_j^{\text{noisy}}(\vec{x}_i) = u_j(\vec{x}_i)  + 0.1 c_{i,j} |\bar{u}_j|, \quad 1\leq i\leq N_{\text{sensors}},\]
where $|\bar{u}_j|$ is the average of the $j$th system's response and the $c_{i,j}\sim \mathcal{N}(0,1)$ are \emph{i.i.d.} The approximation error of the EGF over the domain is $\epsilon = 3.5\%$ for the POD method, and $\epsilon = 2.7\%$ for the randomized SVD method. The associated test error is $\epsilon_{\text{test}}=8.5\%$ for the POD method and $8.9\%$ for the randomized SVD. The error does increase, as we increase the noise in the data set; though, the POD method still learns a reasonable Green's function approximation from the noisy dataset. This suggests that the POD model is not extremely sensitive to noise, in this case. We find that the randomized SVD has similar noise sensitivity to the POD technique, despite the fact that the response are perturbed twice: first when solving the PDE for random forcing terms and secondly for solution under the orthonormal function excitations.

\subsubsection{Effect of different hyperparameters on the POD method} \label{effect-parmeter}

\begin{figure}[htbp]
\centering
\begin{overpic}[width=0.8\textwidth]{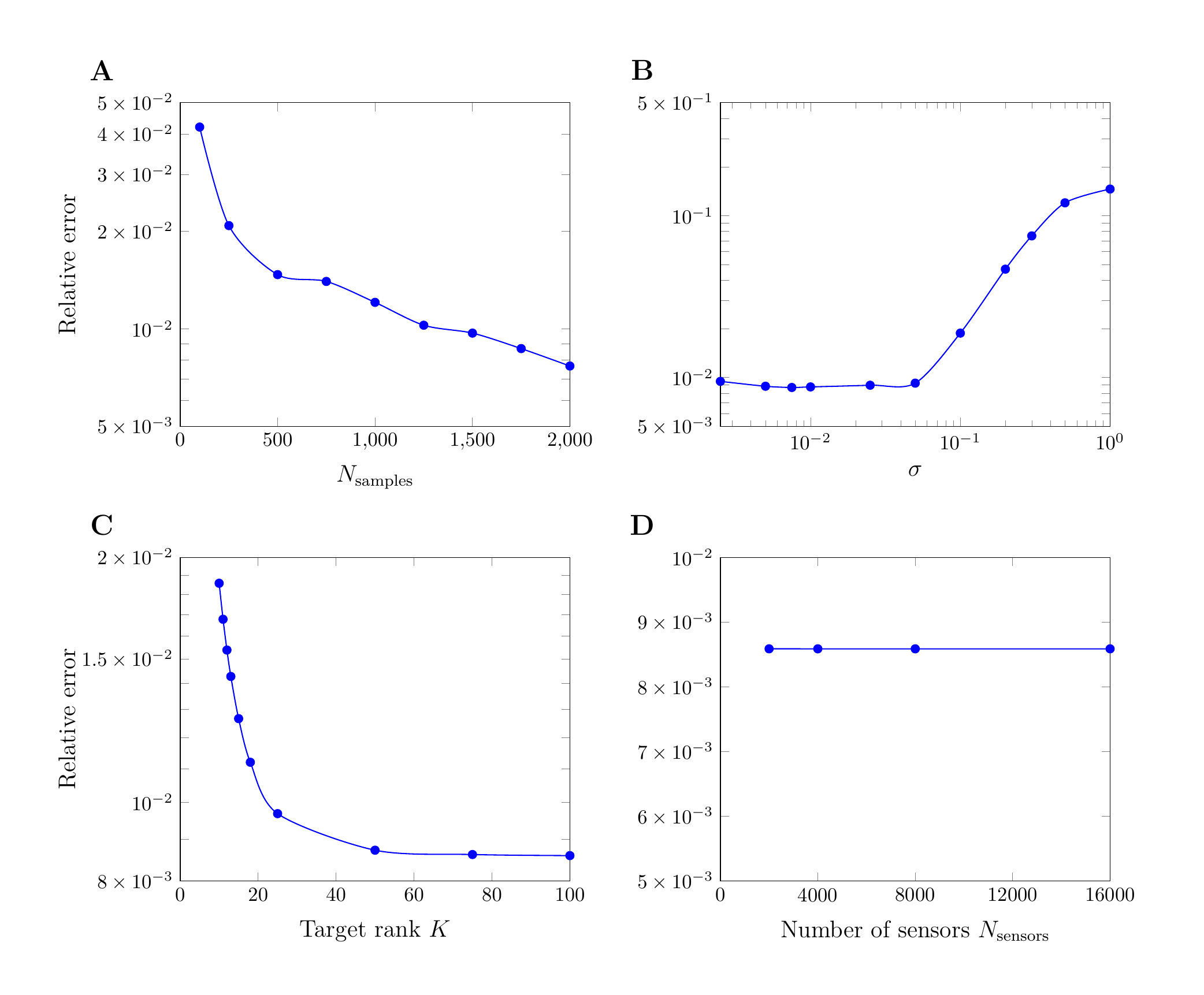}
\end{overpic}
\caption{Relative error for the Green's function associated with the 1D Poisson problem as a function of the following parameters: $N_{\text{samples}}$ (A), $\sigma$ (B), $K$ (C), and $N_\text{sensors}$ (D).}
\label{parameteric-study}
\end{figure}

There are three dominant parameters influencing the fidelity of our empirical Green's functions in the case of the POD method: the number of input-output pairs used for learning the empirical Green's function, the length-scale parameter of the squared-exponential covariance kernel of the Gaussian Process, and finally the number of modes employed within the EGF. \edit{The number of sensors in the domain is chosen so that the forcing functions are sufficiently resolved. We illustrate the effects of these parameters on the EGF's fidelity for the 1D Poisson problem. We report the relative error between the learned and exact Green's function as a function of the parameters in \cref{parameteric-study}. We discretize the domain with quadratic Lagrange ($L_2$) finite elements. When the parameters are not varied for the study, they are fixed at $N_{\text{samples}} = 2000$, $\sigma = 0.0025$, $K = 100$, and $N_{\text{sensors}} = 2000$.} For every set of parameters, we average the results over 10 simulations with different random forcing functions.

\cref{parameteric-study}(A) shows that the error decreases exponentially fast as the number of training pairs in the dataset increases up to $N_{\text{samples}}=500$, where it reaches a slower regime. Additionally, as we decrease the length scale parameter, the forcing functions sampled from the GP become more oscillatory. This allows the forcing functions to more effectively probe the domain of the EGF solution operator, and generate a diverse output ensemble, which increases the quality of the empirical eigenfunctions~\cite{Boulle2021Theory}. Learning empirical eigenvectors with higher fidelity to the eigenfunctions of the exact Green's function leads to more accurate EGF models.

When we plotting the error as a function of the target rank $K$, we find that the error decreases while $K$ increased up to $K=50$. As one allows more POD modes to be used in the representation of the EGF, the EGF accuracy improves up to that point. Since $K$ is effectively the rank at which we are truncating the singular value decomposition, one can choose a value of $K$ based on the decay in singular values. The singular value decay sheds light on where the plateau will form in the panel C of the plot in \cref{parameteric-study}.

\edit{Finally, in \cref{parameteric-study}(D), we plot the relative error as a function of the discretization size. We observe that error stays constant as we increase $N_{\text{sensors}}$ from a level which allows for sufficient resolution of the forcing functions on the domain. This demonstrates that the empirical Green's function which the method learns is independent of the level of discretization of the domain.}

\subsubsection{``Multi-Physics'' context}

\begin{figure}[htbp]
\centering
\begin{overpic}[width=0.8\textwidth]{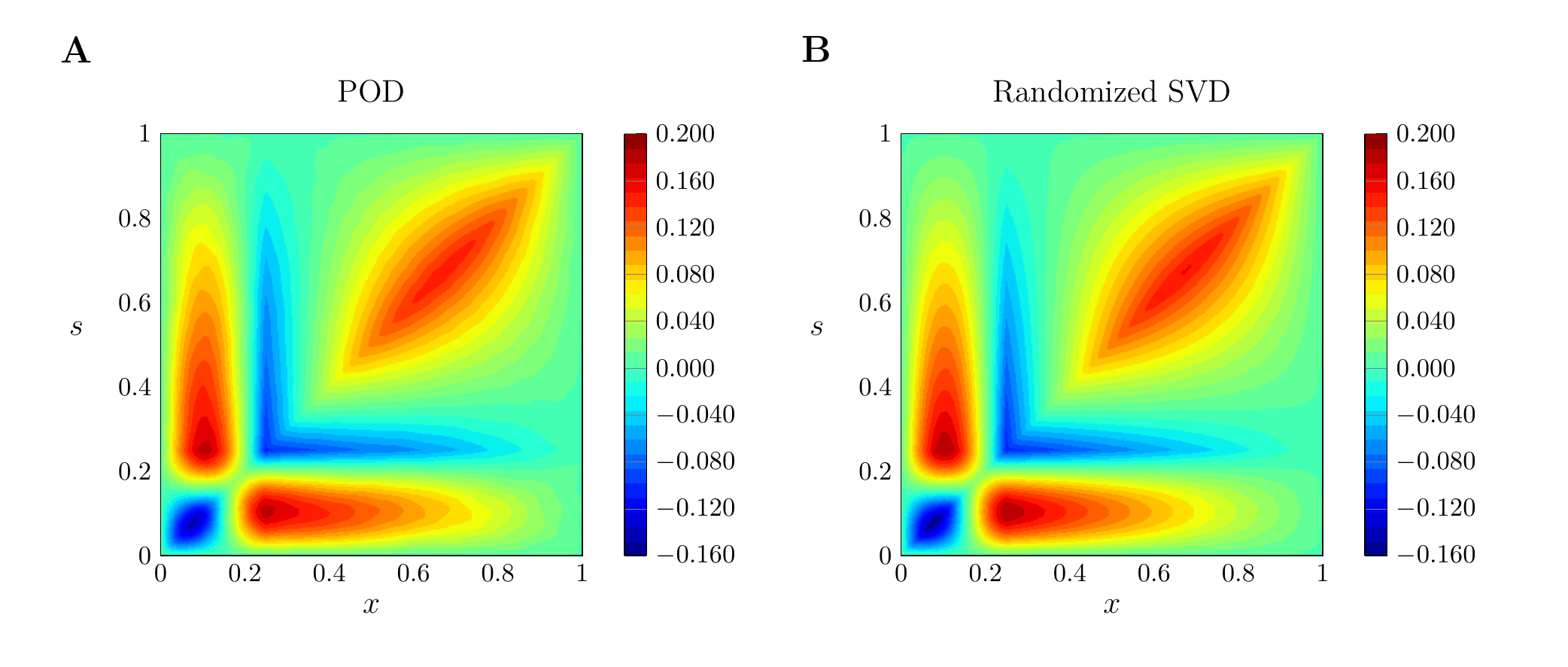}
\end{overpic}
\caption{Green's function associated with \cref{eq_MP} approximated by the POD method (A) and the randomized SVD method.}
\label{multi-physics}
\end{figure}

We now demonstrate our method on a problem in which the behavior of the physical context changes within the domain. The differential equation governing the system on the domain $\Omega = [0,1]$ is given by
\begin{equation} \label{eq_MP}
\begin{gathered}
\frac{1}{2}\left(\frac{d^2 u}{dx^2} + \theta^2 u \right) = f \quad  \text{on}\quad (0,1/4),\qquad -\frac{d^2 u}{dx^2} = f \quad \text{on}\quad (1/4,1),\\
u(0) = u(1/4) = u(1) = 0.
\end{gathered}
\end{equation}
As a result of this, the system is governed by a Helmholtz equation with parameter $\theta = 15$ on the left side of the domain, while it behaves as a Poisson equation for $x>1/4$. We display the learned EGF obtained from the POD method and the randomized SVD method in the left and right panel of \cref{multi-physics}, respectively. The relative error on the testing data set is given by $\epsilon_{\text{test}}=5.8\%$ for POD method and $\epsilon_{\text{test}}=0.3\%$ for the randomized SVD method.

\subsubsection{Summary of errors for 1D problems} \label{sec_error_1d}

The errors for all of the 1D problems are summarized in \cref{summary-errors}.
For all problems, there is a version with clean data and another version in which the system responses are contaminated with $10\%$ additive white Gaussian noise (similar to the 1D Poisson problem in \cref{poisson1D-noisy-section}).

    \begin{table}[htbp]
    \centering
\caption{Summary of the test error with and without noise for the 1D Poisson, Helmholtz, Airy, and ``multi-physics'' problems. The relative error is mentioned in parenthesis for the cases where the exact Green's function is known.}
\vspace{0.3cm}
\begin{tabular}{llcccc}
\toprule
& \multicolumn{1}{c}{Algorithm}  & Poisson & Helmholtz & Airy & Multi-Phys. \\
\midrule
\multicolumn{1}{l}{\multirow{2}{*}{Clean}} & POD & 1.6 (1.1) & 4.0 (1.9) & 2.4 & 5.9 \\
\multicolumn{1}{l}{} & Randomized SVD  & 0.1 (0.09) & 0.2 (0.06) & 0.4 & 0.3 \\
\midrule
\multicolumn{1}{l}{\multirow{2}{*}{Noisy}} & POD & 8.5 (3.5) & 4.6 (1.9) & 6.5 & 8.6 \\
\multicolumn{1}{l}{} & Randomized SVD  & 8.9 (2.7) & 2.4 (0.3) & 7.1 & 1.2 \\
\bottomrule
\end{tabular}
\label{summary-errors}
\end{table}

Looking at \cref{summary-errors}, we observe relatively low test errors (and relative errors), in the noise-free and also noise-contaminated cases. \edit{In general, the randomized SVD achieves better performance than the POD method. One possible explanation for this is that the minimization problem of the POD method (see \cref{minimization-eq}) may exhibit multiple local minima~\cite{Earlsf}. As a result, the POD method might fail to identify the optimal set of singular values associated with the high frequency modes, whereas the randomized SVD does a fairly good job at approximating them, even in the presence of noise. However, there are a few examples (see the Poisson and Airy problems in \cref{summary-errors}), where POD performs marginally better than the randomized SVD when the inputs are contaminated with noise. This could be attributed to the fact that the responses are perturbed twice with the randomized SVD method (during the sketching and projection steps), whereas the responses are only perturbed once in the POD method.}

\subsubsection{Two dimensional Poisson equation} \label{sec_2d_poisson}

\begin{figure}[htbp]
\centering
\vspace{0.2cm}
\begin{overpic}[width = \textwidth]{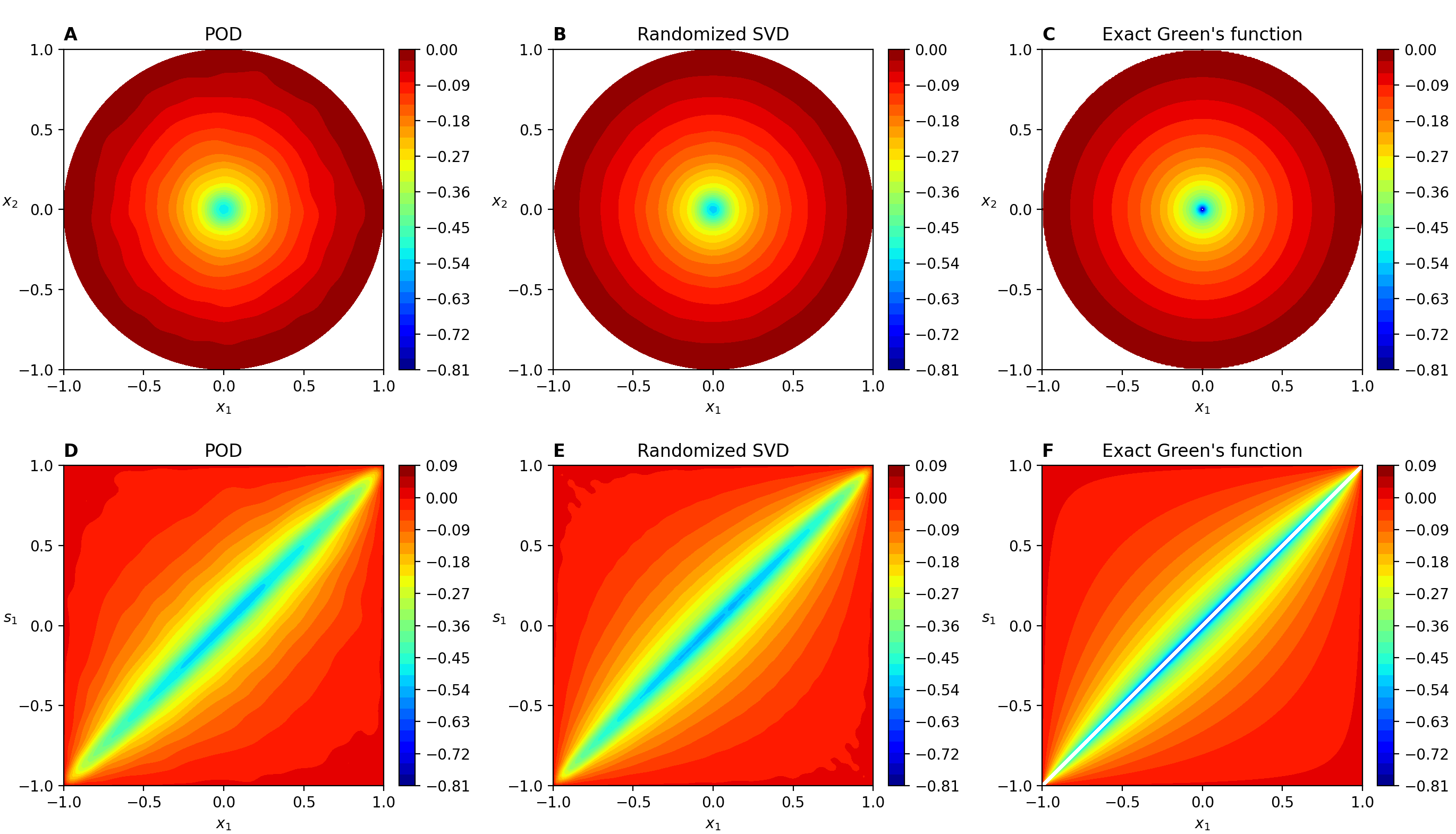}
\end{overpic}
\caption{Green's function associated with the two-dimensional Poisson equation approximated by the POD method (A,D) and randomized SVD method (B,E), along with the exact Green's function (C,F). The top row displays the slice $G(x_1,x_2,0,0)$ while the bottom row shows the Green's function at $x_2=s_2=0$, \emph{i.e.}~$G(x_1,0,s_1,0)$.}
\label{poisson2D}
\end{figure}

Finally, the two algorithms for approximating Green's functions are compared on a Poisson problem in two dimensions. The governing equation with homogeneous boundary conditions on the unit disk reads:
\begin{equation} \label{eq_Poisson}
\begin{aligned}
\nabla \cdot ( \nabla u) &= f, \quad && \text{in } \Omega=D(0,1),\\
u &= 0, \quad &&\text{on } \partial\Omega.
\end{aligned}
\end{equation}
In this case, the closed-form expression for the associated Green's function is known and  given by~\cite{myint2007linear}:
\[G_{\text{exact}}(\vec{x},\vec{s}) = \frac{1}{4 \pi} \ln\left({\frac{(x_1-s_1)^2 + (x_2-s_2)^2}{(x_1 s_2 - x_2 s_1)^2 + (x_1 s_1 + x_2 s_2 - 1)^2}}\right),\quad \vec{x}\neq \vec{s}\in \Omega,\]
where $\vec{x}=(x_1,x_2)$ and $\vec{s}=(s_1,s_2)$. Following \cref{subsec-building-ensemble}, the forcing functions $\{f_j\}_{j=1}^{N_{\text{samples}}}$ are generated using the \texttt{randnfundisk} function of the Chebfun software system, with a length-scale parameter of $\sigma = 0.2$. \cref{eq_Poisson} is discretized with quadratic Lagrange finite elements on a mesh of $N_{\text{sensors}}\approx 5000$ nodes and solved with the FEniCS software. We select a target rank of $K=500$ and choose $N_{\text{samples}}=2000$ for the POD method, while the randomized SVD method uses $N_{\text{samples}}=500$ input-output pairs. After computing an approximation to the Green's function using both methods, we visualize the two-dimensional slices at $s_1=s_2=0$ and $x_2=s_2=0$ of the learned EGF and exact Green's function in \cref{poisson2D}. We subsequently observe a very good agreement between the low-rank EGF and the exact Green's function, as confirmed by the low relative test errors of $\epsilon_{\text{test}} = 4.7\%$ for the POD method and $\epsilon_{\text{test}} = 0.4\%$ for the randomized SVD technique. We remark that the presence of a pole at $x=s$ leads to a slow decay in the singular values associated with the Green's function, since the function is not smooth and requires a larger number of training pairs to obtain an accurate approximation.

\subsection{Interpolation and extrapolation of Green's functions}

We now evaluate the interpolation algorithm described in \cref{sec-interpolation}.

\subsubsection{One dimensional Airy problem} \label{subsec-airy-interp}

\begin{figure}[htbp]
\centering
\begin{overpic}[width=\textwidth]{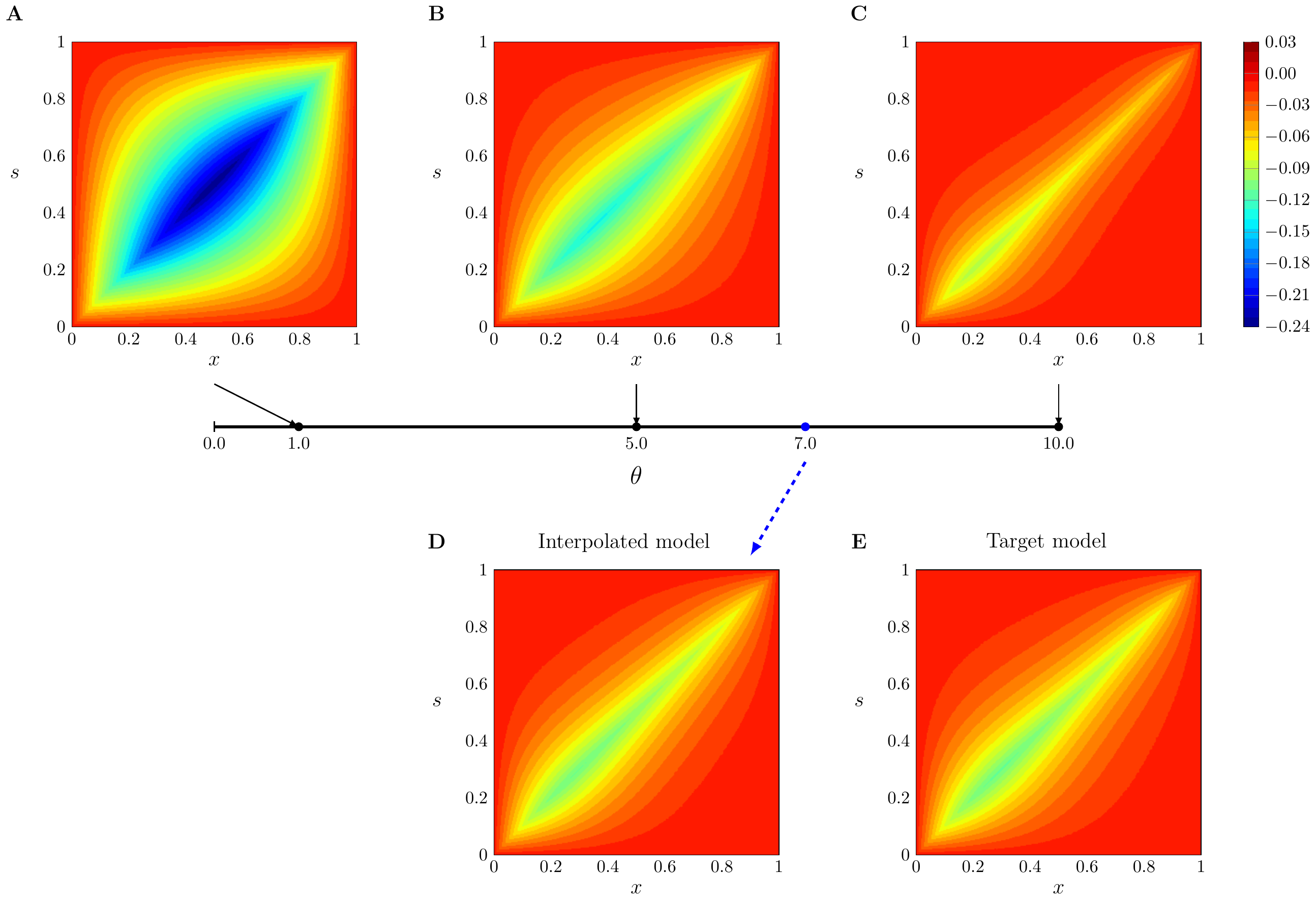}
\end{overpic}
\caption{(A-C) Approximation of the Green's functions associated with the Airy problem at $\theta=1,5,10$ used by the interpolation scheme. (D) Interpolated Green's function at $\theta_*=7$. (E) The target Green's function approximated by randomized SVD.}
\label{interp-airy}
\end{figure}

We begin with a parameterized Airy's equation in one spatial dimension with homogeneous Dirichlet boundary conditions:
\begin{equation} \label{eq_airy}
\frac{d^2u}{dx^2} - \theta^{2} x u = f, \quad x \in [0,1],
\end{equation}
where $\theta\in \R$ is a parameter of the model. We first compute approximations to the Green's functions at parameter values $\theta_1 = 1$, $\theta_2 = 5$, and $\theta_3 = 10$, using the randomized SVD. Then, we aim to interpolate the Green's function at $\theta_* = 7$, where there is no training data. Since we do not have access to an analytical expression for the Green's function associated with \cref{eq_airy}, we compare the interpolated model against the approximation (target Green's function) computed by the randomized SVD at the target parameter in \cref{interp-airy}, \emph{i.e.}, using a dataset generated at $\theta=\theta_*$. We obtain a relative error of $\epsilon=2.7\%$ between the target and interpolated Green's functions, while the relative error on the testing data set at $\theta_*=7$ is equal to $\epsilon_{\text{test}}=2.6\%$. While the error is six times larger than the one reported in \cref{sec_error_1d} for the randomized SVD method, we note that the error is on the order of what the POD method reached. This is a relatively low error given the very small number of interpolation points used, \emph{i.e.}, just three points.

\subsubsection{Extrapolation on 1D Airy problem}

\begin{figure}[htbp]
\centering
\begin{overpic}[width=0.8\textwidth]{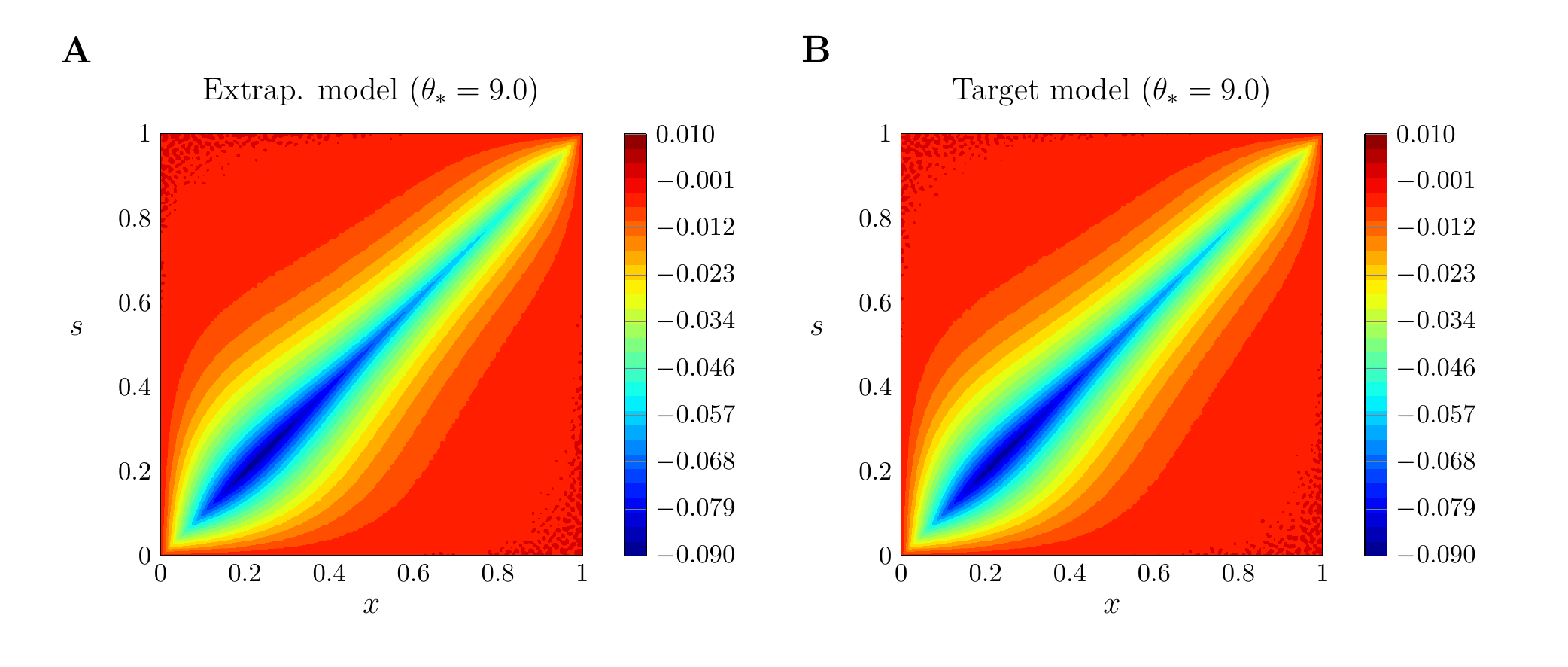}
\end{overpic}
\caption{EGFs for 1D Airy problem at $\theta_* = 9.0$ obtained using extrapolation on a manifold (A) and learned from data at the target parameter (B).}
\label{extrapolation-comparison}
\end{figure}

We note that our method also generalizes to extrapolation, out of the neighborhood of a set of learned empirical eigenmodes. To demonstrate this we compute approximations to the Green's functions at parameter values of $\theta_1 = 6.0$, $\theta_2 = 7.0$, and $\theta_3 = 8.0$ using the randomized SVD, and then subsequently extrapolate the EGF model to $\theta_* = 9.0$. We compare the extrapolated EGF against a randomized SVD based EGF model, learned (using ground truth, system response data) at the target parameter $\theta_*$ in \cref{extrapolation-comparison}.

We obtain a relative error of $\epsilon = 1.3\%$ between the target and extrapolated EGFs and testing error for the extrapolated EGF of $\epsilon_{\text{test}} = 1.1\%$. These errors are similar to the interpolation problem. While extrapolation is in general less reliable than interpolation, our method does appear to perform reasonably well in this case.

\subsubsection{2D Helmholtz problem}

\begin{figure}[htbp]
\centering
\begin{overpic}[width = 0.8\textwidth]{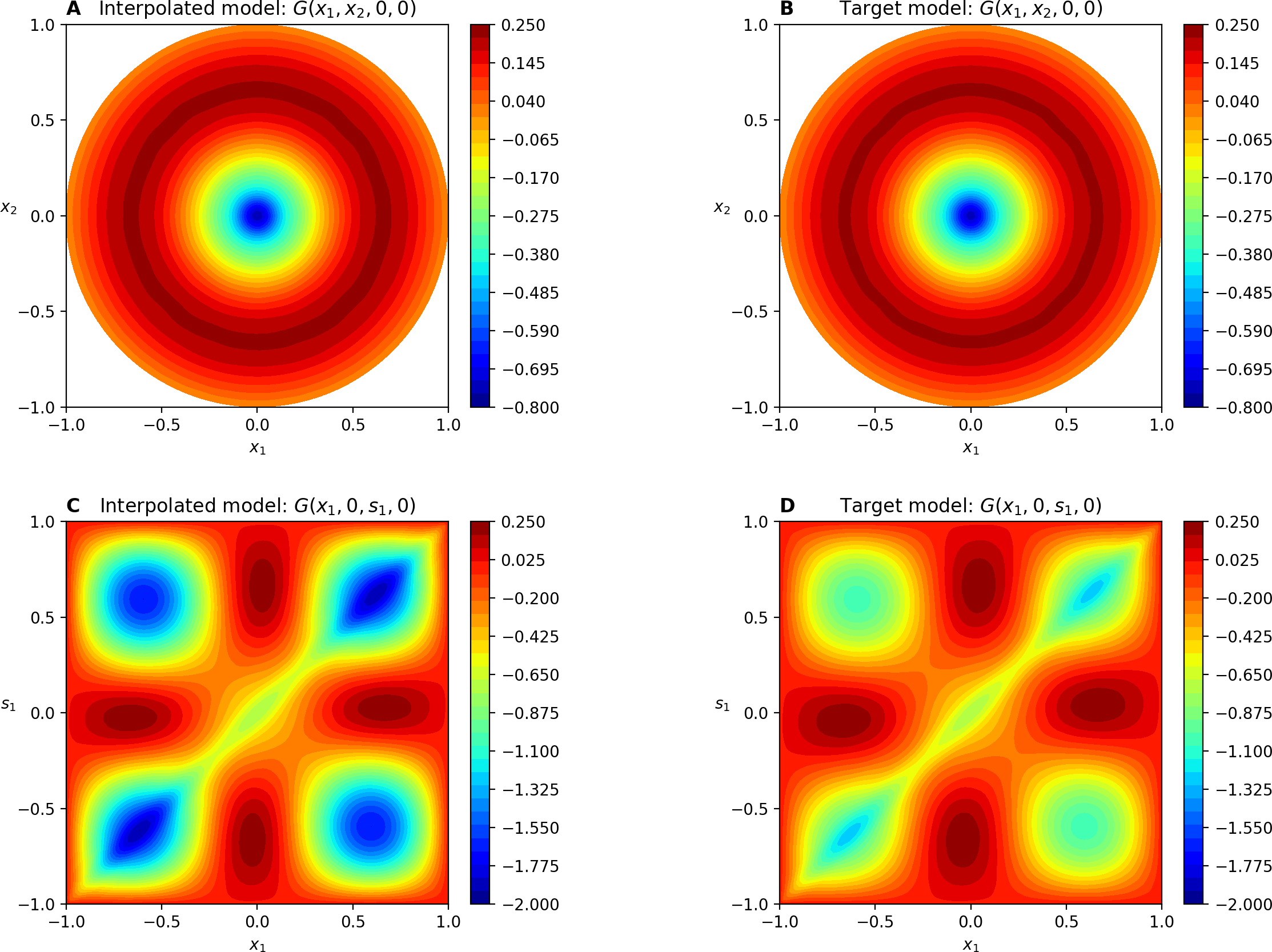}
\end{overpic}
\caption{Green's function associated with the two-dimensional Helmholtz problem at $\theta_* = 5.0$, obtained by manifold interpolation (A,C), and learned at the target parameter from data (B,D). The top row display the slice $G(x_1,x_2,0,0)$, while the bottom row shows the slice $G(x_1,0,s_1,0)$.
}
\label{interp-helm2D}
\end{figure}

\begin{figure}[htbp]
\centering
\begin{overpic}[width = \textwidth]{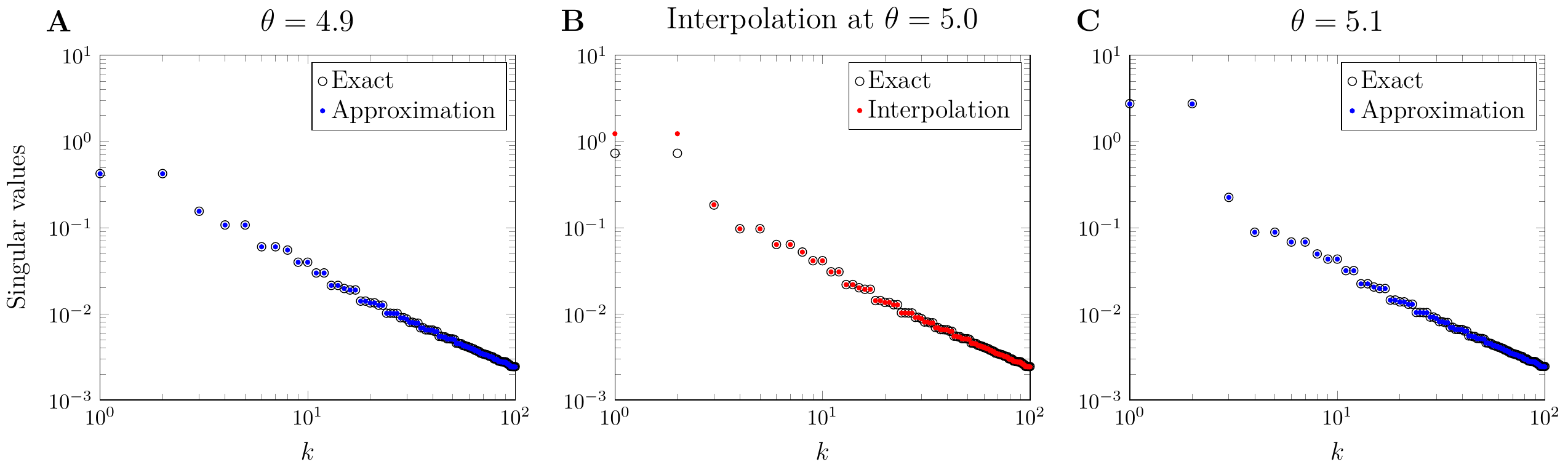}
\end{overpic}
\caption{First hundred largest singular values of the Green's function with the two-dimensional Helmholtz problem approximated by the randomized SVD method at $\theta=4.9$ (A) and $\theta=5.1$ (C). Panel (B) displays the first hundred interpolated singular values along with the exact ones at $\theta_*=5.0$.}
\label{interp-svd_h_2d}
\end{figure}

We now demonstrate our interpolation method on a 2D Helmholtz problem. The governing equations are given by
\begin{equation} \label{eq_2D_helm}
    \begin{aligned}
        \nabla \cdot ( \nabla u) + \theta^2 u &= f, \quad && \text{in } \Omega=D(0,1),\\
        u &= 0, \quad &&\text{on } \partial\Omega.
    \end{aligned}
\end{equation}
We employ the randomized SVD with a length-scale parameter of $\sigma = 0.2$, $N_{\text{samples}} = 100$, and the number of empirical eigenmodes used in the EGF model fixed to $K = 100$. We solve \cref{eq_2D_helm} with piecewise quadratic finite elements using FEniCS. The nodes of the mesh define the sensor locations $\bar{\bm{x}} \in \Omega$. The number of sensors is then, $N_{\text{sensors}} \sim 10000$. We learn empirical eigenmodes at parameter values of $\theta_1 = 4.8$, $\theta_2 = 4.9$, and $\theta_2 = 5.1$. Using these empirical eigenmodes, we subsequently apply our interpolation method and obtain an EGF at $\theta_* = 5.0$. A comparison between the interpolated Green's function (left) and a Green's function learned using a randomized SVD based EGF model (\emph{i.e.}, learned using ground truth system response data at the target parameter $\theta_*$) is shown in \cref{interp-helm2D}. For this problem, the largest singular values change rapidly near $\theta_{\text{crit}} \approx 5.14$. While we are only using three points to find an interpolated Green's function at $\theta_* = 5.0$, resulting in relatively high errors for the largest interpolated singular values, the remaining singular values are interpolated with high accuracy, as shown in \cref{interp-svd_h_2d}.

\edit{
\subsubsection{One dimensional fractional laplacian problem} \label{subsec-laplacian-interp}

In this section, we evaluate the interpolation method on a one-dimensional fractional Laplacian operator with periodic boundary conditions:
\begin{equation}\label{eq_fraclaplacian}
(-\Delta)^{\theta} u = f, \quad u(-1)=u(1),\quad x \in [-1,1], \\
\end{equation}
where $0<\theta<1$ is the fractional order. We solve \cref{eq_fraclaplacian} using a Fourier spectral collocation method and use the randomized SVD to  approximate the Green's function at parameters $\theta_1 = 0.6$, $\theta_2 = 0.7$, and $\theta_3 = 0.8$. Then, we interpolate the Green's functions at a new parameter $\theta_* = 0.75$ for which there is no training data available. We compare the interpolated model against the approximant (target Green's function) computed by the randomized SVD at the target parameter $\theta_*$ in \cref{interp-fracplaplacian}, \emph{i.e.}, using a dataset generated at the target parameter. We obtain a relative error of $\epsilon= 0.1\%$ between the target and interpolated Green's function, while the relative error on the testing dataset at $\theta_*=0.75$ is equal to $\epsilon_{\text{test}}=0.8\%$.

\begin{figure}[htbp]
\centering
\begin{overpic}[width=\textwidth]{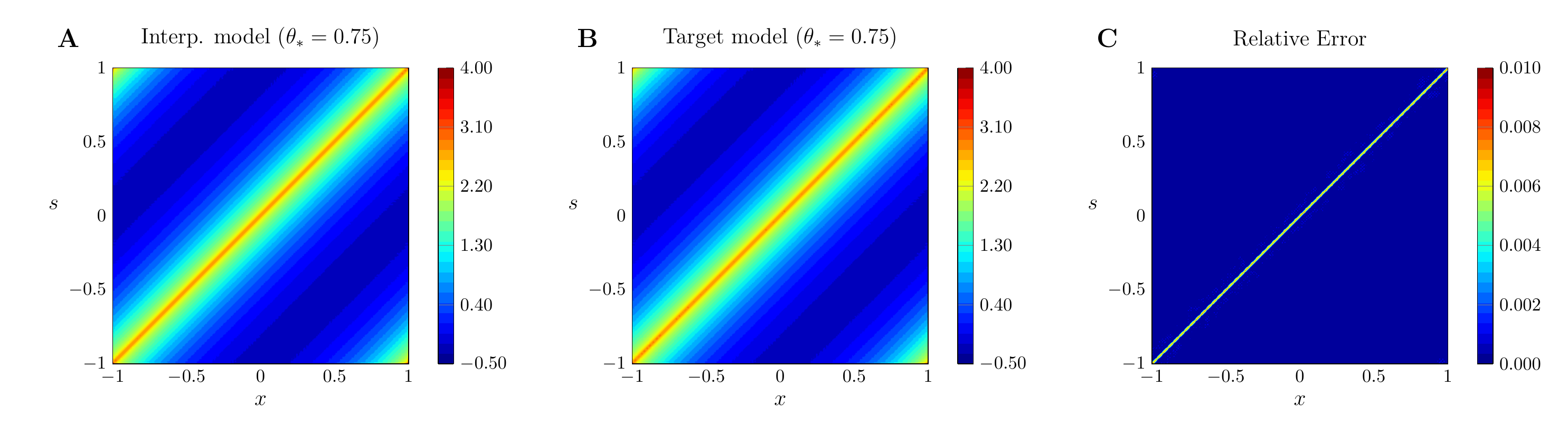}
\end{overpic}
\caption{(A) Approximation of the Green's function associated with the fractional Laplacian problem at $\theta_* = 0.75$ used by the interpolation scheme. (B) Interpolated Green's function at $\theta_*=0.75$. (C) Relative error between the interpolated and target model.}
\label{interp-fracplaplacian}
\end{figure}

\subsubsection{2D Poisson on a Pac-Man-shaped domain}

In this last example, we consider a more challenging parametrization, where the shape of the domain depends on the parameter $\theta$. We consider a two-dimensional Poisson problem defined on a Pac-Man-shaped domain $P(\theta)\subset \mathbb{R}^2$, with a variable-size inscribed angle $\theta\in [0,2\pi)$. The governing equations are given by:
\begin{equation} \label{eq_pacman}
    \begin{aligned}
        \nabla \cdot ( \nabla u) &= f, \quad && \text{in } \Omega = P(\theta),\\
        u &= 0, \quad &&\text{on } \partial\Omega.
    \end{aligned}
\end{equation}
We approximate Green's functions associated with \cref{eq_pacman} using the randomized SVD at parameters $\theta_1 = 0.4 \pi$, $\theta_2 = 0.5 \pi$, and $\theta_3 = 0.6 \pi$, and apply our interpolation algorithm to recover the Green's function at $\theta_*=0.55\pi$. To ensure that the empirical eigenmodes are defined on the same finite element basis, we first generate a mesh of $P(\theta)$ at $\theta_*$ using the Gmsh meshing software~\cite{gmsh}. Then, we apply a polar transformation of the coordinates of the mesh to generate meshes at $\theta_1$, $\theta_2$, and $\theta_3$ (see \cref{mesh-pacman}). If $(r_*, \alpha_*)$ represents the polar coordinates associated the mesh $P(\theta_*)$ and $(r_1, \alpha_1)$ represents the same for the mesh $P(\theta_1)$, then the coordinate transform for changing a mesh $P(\theta_*)$ to $P(\theta_2)$ is defined as
\[
r_1 = r_*, \quad \alpha_1 =
\begin{cases}
a \alpha_* + b, & \text{if}\ r \sin(\alpha_*) \ge 0, \\
a \alpha_* + b, & \text{if}\ r \sin(\alpha_*) < 0,
\end{cases}
\]
where $a = \frac{\pi - \alpha_1/2}{\pi - \alpha_*/2}$ and $b = \pi - a \pi$. We solve \cref{eq_pacman} using the Firedrake finite element software~\cite{Rathgeber2016} and employ the randomized SVD with the parameters $\sigma = 0.2$, $N_{\text{samples}} = 200$, $N_{\text{sensors}} \sim 5000$, $K = 200$, to approximate the Green's functions at parameters $\theta_1$, $\theta_2$, and $\theta_3$.

\begin{figure}[htbp]
\centering
\begin{overpic}[width = \textwidth]{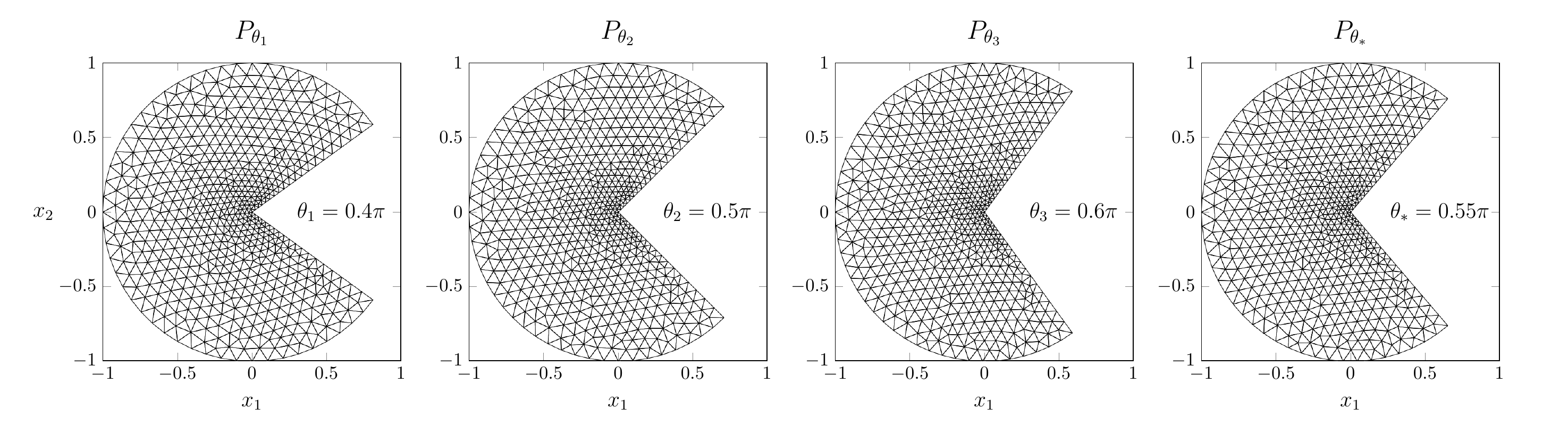}
\end{overpic}
\caption{Meshes of the Pac-man shaped domain $P(\theta)$ with an inscribed angle $\theta$.}
\label{mesh-pacman}
\end{figure}

In this case, the eigenmodes for the empirical Green's functions learned at different parameters are defined on different domains. We then map the different eigenmodes to the target domain $P(\theta_*)$ using the coordinate transform defined above. Here, we exploit Firedrake's capabilities to compute the Jacobian associated with the coordinate transformation. We multiply the eigenmodes ($\bm{\Psi}(\bar{\bm{x}};\theta_j)$) for the interpolant with the determinant of the inverse of the Jacobian, to preserve the orthonormality of modes with respect to the mesh at $\theta_*$. Using these empirical eigenmode sets, we subsequently apply our interpolation method and obtain an EGF at $\theta_* = 0.55 \pi$.

\begin{figure}[htbp]
\centering
\begin{overpic}[width = \textwidth]{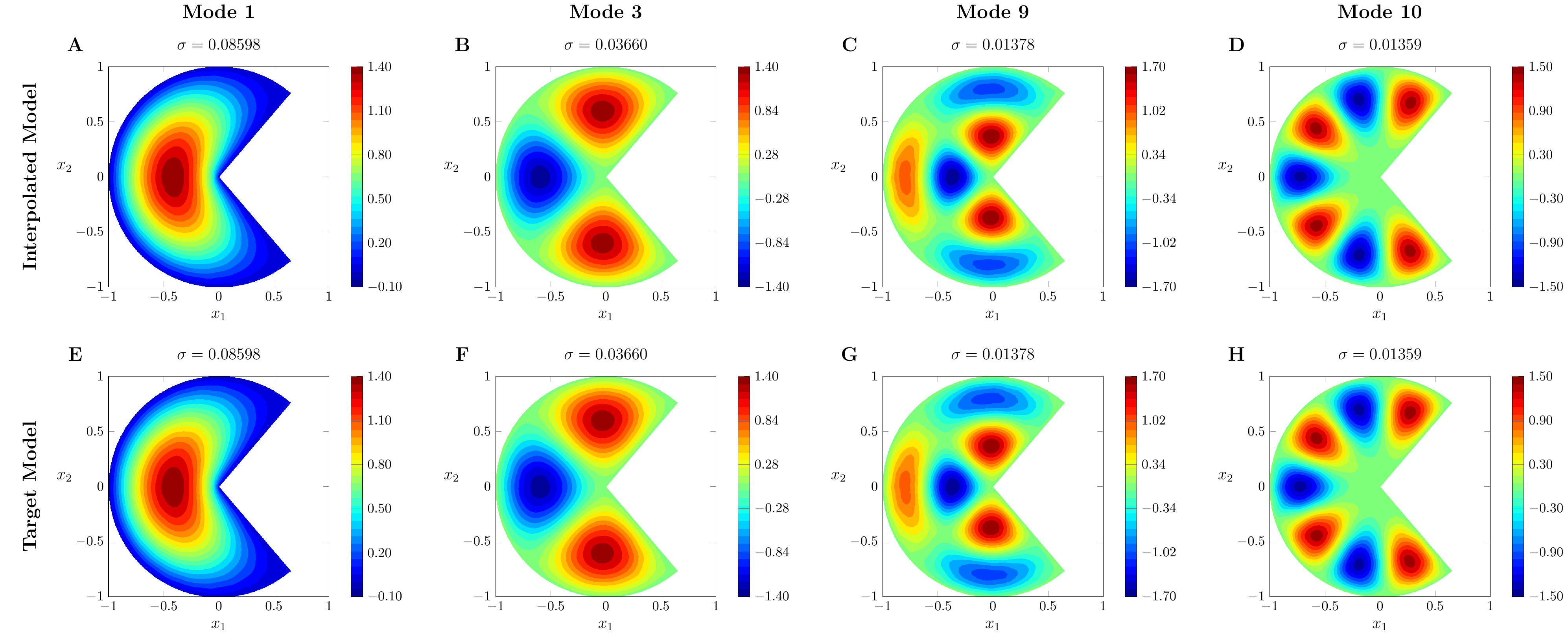}
\end{overpic}
\caption{Comparison between the modes of the Green's function associated with the two-dimensional Poisson problem on a parameterized Pac-man shaped domain, approximated by interpolated empirical Green's functions at $\theta_*=0.55\pi$ (A,B,C,D) and the randomized SVD (E,F,G,H).}
\label{interp-pacman}
\end{figure}

\begin{figure}[htbp]
\centering
\begin{overpic}[width = 0.5\textwidth]{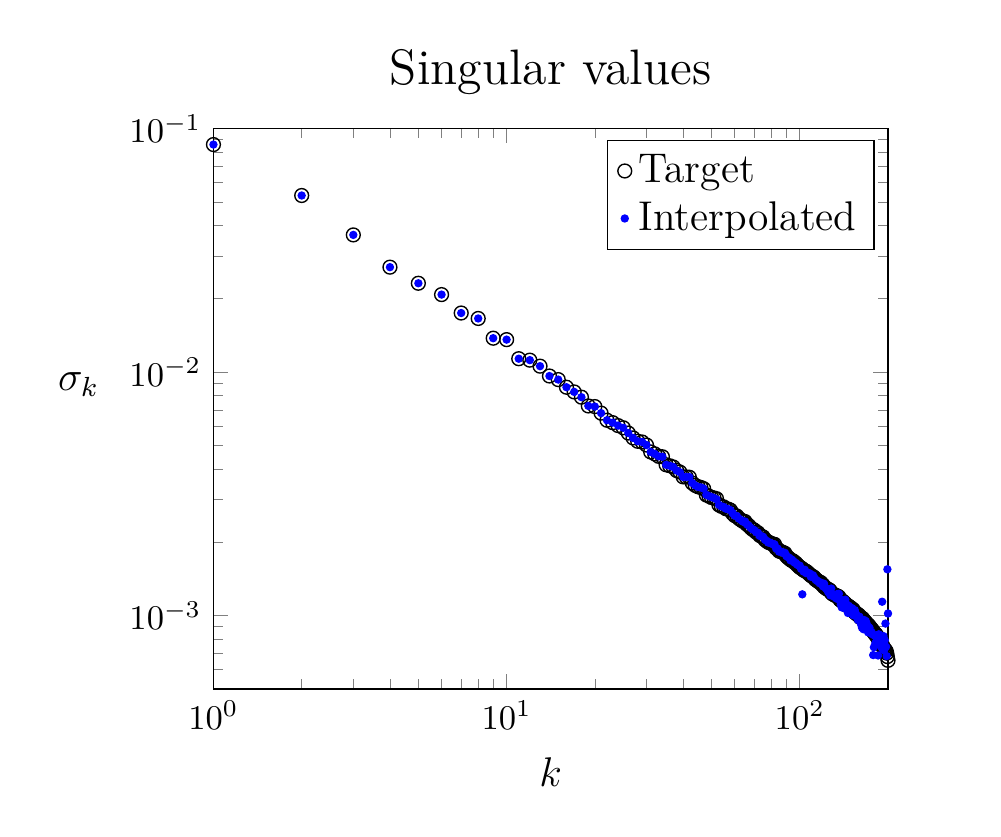}
\end{overpic}
\caption{Two hundred largest singular values of the interpolated Green's function associated with the two-dimensional Poisson problem on a parameterized Pac-man shaped domain, compared against a target Green's function approximated using the randomized SVD at $\theta_*=0.55\pi$.}
\label{interp-svalues-pacman}
\end{figure}

A comparison between the modes of the interpolated Green's function (first row) and an approximation (target Green's function) computed by the randomized SVD method (second row) at the target parameter (\emph{i.e.},~using a dataset generated at $\theta_*$), is shown in \cref{interp-pacman}. A comparison between the coefficients learned using a randomized SVD based empirical Green's function (\emph{i.e.},~learned using ground truth, system response data at the target parameter) and the interpolated empirical Green's function is shown in \cref{interp-svalues-pacman}. Note that we reorder the eigenmodes for the interpolated Green's function to match the eigenmodes of the target empirical Green's function, as described in \cref{sec_correct_sign_order}.
}

\section{Conclusions} \label{sec-conclusion}

In this work, we considered two methods based on the proper orthogonal decomposition and the randomized SVD for learning empirical Green's functions from training data consisting of excitation-response pairs. The first method is more appropriate in cases where one has little control over the forcing functions, while the latter is more accurate, but requires more control over the forcing terms and two passes through the PDE. Both methods are observed to perform well in the one and two dimensional numerical experiments. Then, we proposed the use of a manifold interpolation scheme in an offline-online setting, where offline excitation-response data, taken at specific model parameter instances, are compressed into empirical eigenmodes. These eigenmodes are subsequently used within a manifold interpolation scheme, to uncover other suitable eigenmodes, for an unseen model parameter instance; thus rendering an online, ``just-in-time'' EGF (obtained without the benefit of excitation-response data) This interpolation approach is demonstrated in 1D and 2D contexts; yielding promising results.

\section*{Data availability}

The code used to produce the numerical results is publicly available on GitHub at \url{https://github.com/hsharsh/EmpiricalGreensFunctions} for reproducibility purposes.

\section*{Acknowledgements}

H.P. and C.J.E. were supported by the Army Research Office (ARO) Biomathematics Program grant W911NF-18-1-0351. N.B. was supported by the EPSRC Centre for Doctoral Training in Industrially Focused Mathematical Modelling through grant EP/L015803/1 in collaboration with Simula Research Laboratory and an INI-Simons Postdoctoral Research Fellowship. The authors also thank Max Jenquin for initial help in implementing this idea, as well as Maria Oprea, for insights regarding the manifold interpolation. We are grateful to Alex Townsend for his insights and many helpful suggestions.

\bibliographystyle{elsarticle-num}
\bibliography{citations}
\end{document}